\documentclass[english,11pt]{amsart}
\usepackage{sty/preamble} 

\newcommand{\papertitle}{May's Conjecture on Bimonoidal Functors and Multiplicative Infinite Loop Space Theory}

\title{\papertitle}
\author{Donald Yau}
\address{The Ohio State University at Newark}
\email{yau.22@osu.edu}
\date{\today}


\begin{document}
\thispagestyle{empty}

\begin{abstract}
A conjecture of May states that there is an up-to-adjunction strictification of symmetric bimonoidal functors between bipermutative categories.  The main result of this paper proves a weaker form of May's conjecture that starts with multiplicatively strong symmetric bimonoidal functors.  As the main application, for May's multiplicative infinite loop space machine from bipermutative categories to either $\Einf$-ring spaces or $\Einf$-ring spectra, multiplicatively strong symmetric bimonoidal functors can be replaced by strict symmetric bimonoidal functors.
\end{abstract}

\maketitle
\tableofcontents

\section{Introduction}\label{sec:introduction}

The main result of this paper proves an up-to-adjunction strictification of multiplicatively strong symmetric bimonoidal functors between bipermutative categories, confirming a weaker form of a conjecture of May.  The categorical context and application to multiplicative infinite loop space theory are explained after the following brief review of infinite loop space theory.

\subsection*{Permutative Categories and Infinite Loop Space Machine}

May's infinite loop space machine $\mayk$ \cite{may,may-groupcompletion,may-permutative} is a functor that sends small symmetric monoidal categories to $K$-theory connective spectra.  By Isbell's Strictification Theorem of symmetric monoidal categories \cite{isbell}, there is no loss of generality in restricting attention to small permutative categories.  At the morphism level, May's machine $\mayk$ applies to symmetric monoidal functors.  Analogous to the input data at the object level, it would be desirable if symmetric monoidal functors can be replaced by \emph{strict} symmetric monoidal functors.  This is possible by the following theorem from \cite[Prop.\ 4.3]{may-pairings}, which proves an up-to-adjunction strictification of symmetric monoidal functors between permutative categories. 

\begin{theorem}[May]\label{may4.3}
There is a functor
\[\perm \fto{\Pst} \permst\]
from the category $\perm$ of small permutative categories and symmetric monoidal functors to the category $\permst$ of small permutative categories and strict symmetric monoidal functors.  Moreover, for each permutative category $\A$, there is a monoidal adjunction
\[\begin{tikzpicture}[vcenter]
\def\b{22} \def\h{1.8}
\draw[0cell]
(0,0) node (a1) {\phantom{\A}}
(a1)++(-.1,0) node (a1') {\Pst\A}
(a1)++(\h,0) node (a2) {\A}
(a1)++(\h/2,0) node (a0) {\bot}
;
\draw[1cell=.9]
(a1) edge[bend left=\b] node {\lefta} (a2)
(a2) edge[bend left=\b] node {\righta} (a1)
;
\end{tikzpicture}\]
such that the following three statements hold.
\begin{enumerate}
\item\label{may4.3_i} The counit of the adjunction, $\lefta\righta \to 1_{\A}$, is the identity monoidal natural transformation.
\item\label{may4.3_ii} The left adjoint $\lefta$ is a strict symmetric monoidal functor and is natural in strict symmetric monoidal functors.
\item\label{may4.3_iii} The right adjoint $\righta$ is a symmetric monoidal functor and is natural in symmetric monoidal functors.  
\end{enumerate}
\end{theorem}

In May's \cref{may4.3}, the unit $1_{\Pst\A} \to \righta\lefta$ of the adjunction $(\lefta,\righta)$ is not an isomorphism, so $(\lefta,\righta)$ is not an adjoint equivalence of categories.  Instead, since $\righta$ admits a left adjoint, its image $\mayk\righta$ under May's machine is a weak equivalence of spectra.  Given a symmetric monoidal functor $\sff \cn \A \to \B$ between small permutative categories, the naturality of the right adjoint $\righta$ yields the left commutative diagram of symmetric monoidal functors below.
\begin{equation}\label{kr_weq}
\begin{tikzpicture}[vcenter]
\def\v{-1.3} \def\h{1.8} \def\g{1} \def\l{1.5}
\draw[0cell]
(0,0) node (a11) {\Pst\A}
(a11)++(\h,0) node (a12) {\A}
(a11)++(0,\v) node (a21) {\Pst\B}
(a12)++(0,\v) node (a22) {\B}
;
\draw[1cell=.9]
(a12) edge node[swap] {\righta} (a11)
(a11) edge node[swap] {\Pst f} (a21)
(a12) edge node {\sff} (a22)
(a22) edge node[swap] {\righta} (a21)
;
\begin{scope}[shift={(\h+\g,\v/2)}]
\draw[0cell]
(0,0) node (a1) {\phantom{\A}}
(a1)++(\l,0) node (a2) {\phantom{\A}}
;
\draw[1cell]
(a1) edge[|->] node {\mayk} (a2)
;
\end{scope}
\begin{scope}[shift={(\h+\l+2*\g,0)}]
\draw[0cell]
(0,0) node (a11) {\mayk\Pst\A}
(a11)++(1.2*\h,0) node (a12) {\mayk\A}
(a11)++(0,\v) node (a21) {\mayk\Pst\B}
(a12)++(0,\v) node (a22) {\mayk\B}
;
\draw[1cell=.9]
(a12) edge node[swap] {\mayk\righta} node {\sim} (a11)
(a11) edge[transform canvas={xshift=1ex}] node[swap] {\mayk\Pst f} (a21)
(a12) edge node {\mayk\sff} (a22)
(a22) edge node[swap] {\mayk\righta} node {\sim} (a21)
;
\end{scope}
\end{tikzpicture}
\end{equation}
Applying May's infinite loop space machine $\mayk$ yields the right commutative diagram of spectra such that each morphism $\mayk\righta$ is a weak equivalence.  Thus, from the point of view of infinite loop space theory, the symmetric monoidal functor $\sff$ can be replaced by the \emph{strict} symmetric monoidal functor $\Pst f \cn \Pst\A \to \Pst\B$.

\subsection*{Bipermutative Categories and Multiplicative Infinite Loop Space Machine}

A symmetric bimonoidal category is equipped with an additive symmetric monoidal structure $\bplus$ and a multiplicative symmetric monoidal structure $\btimes$ that distributes over the additive one in a categorical sense.  For example, the category of modules over a commutative ring is a symmetric bimonoidal category with the direct sum and tensor product of modules.  Laplaza's Coherence Theorem \cite{laplaza,laplaza2} for symmetric bimonoidal categories is an analogue of Mac Lane's Coherence Theorem for symmetric monoidal categories.  It leads to a strictification theorem \cite[Prop.\! 6.3.5]{may-einfinity} that strictifies symmetric bimonoidal categories with invertible distributivity morphisms to bipermutative categories, which have two interacting permutative structures.  

A symmetric bimonoidal functor $\sff \cn \A \to \B$ between symmetric bimonoidal categories is a functor equipped with two symmetric monoidal functor structures, an additive one $\fplus = (\sff,\ftwoplus,\fzeroplus)$ and a multiplicative one $\ftimes = (\sff,\ftwotimes,\fzerotimes)$, that are compatible with the distributivity structures. 
A symmetric bimonoidal functor $\sff$ is \emph{multiplicatively strong} if the multiplicative monoidal and unit constraints, 
\[\sff a \btimes \sff b \fto{\ftwotimes} \sff(a \btimes b) \andspace \bunit \fto{\fzerotimes} \sff\bunit,\] 
are isomorphisms.  It is \emph{strict} if both $\fplus$ and $\ftimes$ are strict, meaning that $\ftwoplus$, $\fzeroplus$, $\ftwotimes$, and $\fzerotimes$ are identities.

May's infinite loop space machine has a multiplicative enhancement \cite{may-multiplicative,may-construction,may-precisely} that sends small bipermutative categories to $\Einf$-ring spaces and then to $\Einf$-ring spectra.  At the morphism level, May's multiplicative machine takes symmetric bimonoidal functors as input.  With applications to multiplicative infinite loop space theory in mind, the problem of extending \cref{may4.3} to the context of bipermutative categories is mentioned in \cite[page 16]{may-multiplicative} and restated as the following conjecture in \cite[13.1]{may-construction}.

\begin{conjecture}[May]\label{may-conjecture}
There is a functor on the bipermutative category level that replaces symmetric bimonoidal functors by strict symmetric bimonoidal functors, in a sense analogous to \cref{may4.3} for permutative categories.
\end{conjecture}

The main result of this paper \pcref{thm:sbfstrict} proves May's \cref{may-conjecture} for symmetric bimonoidal functors that are multiplicatively strong. 

\begin{theorem}\label{thm:main}
There is an up-to-adjunction strictification of multiplicatively strong symmetric bimonoidal functors between bipermutative categories. 
\end{theorem}

A precise version of \cref{thm:main} is \cref{thm:sbfstrict}, which is a bipermutative analogue of May's \cref{may4.3}.  The multiplicatively strong hypothesis in \cref{thm:main} is used in the arrow $\ftwotimesinv$ in the diagram \cref{Bsf_mor_ii} to ensure that the strictified symmetric bimonoidal functor is well defined on morphisms.  It is not known whether there is an improvement of \cref{thm:main} where the multiplicatively strong hypothesis is removed.  

The main application of \cref{thm:main} ($=$ \cref{thm:sbfstrict}) is in multiplicative infinite loop space theory.  Using \cref{thm:main}, the paragraph containing \cref{kr_weq} above admits a bipermutative analogue.  Thus, using May's multiplicative infinite loop space machine, multiplicatively strong symmetric bimonoidal functors between small bipermutative categories can be replaced by \emph{strict} symmetric bimonoidal functors.

\subsection*{Organization and Reading Suggestion}
\cref{sec:bipermcat} reviews (bi)permutative categories, symmetric (bi)monoidal functors, (bi)monoidal natural transformations, and Laplaza's Coherence Theorem.  \cref{sec:strict_sbf} proves \cref{thm:sbfstrict}, which is the precise version of \cref{thm:main} and extends May's \cref{may4.3} to the bipermutative setting.  \cref{rk:Bs_adjoint} discusses some conceptual differences between \cref{may4.3,thm:sbfstrict}.  \cref{sec:applications} applies \cref{thm:sbfstrict} to obtain a bipermutative analogue of \cref{kr_weq} involving May's multiplicative infinite loop space machine.  \cref{rk:Kem} discusses the reason why \cref{thm:sbfstrict} does \emph{not} directly apply to Elmendorf-Mandell multifunctorial $K$-theory.

The application to multiplicative infinite loop space theory in \cref{sec:applications} uses only the statement of \cref{thm:sbfstrict}.  For a quick overview of the main result and application, the reader can first read the definitions about bipermutative categories \pcref{def:bipermcat,def:sbfunctor,def:bimonnat}, followed by the statement of \cref{thm:sbfstrict} and \cref{sec:applications}.

\section{Bipermutative Categories}\label{sec:bipermcat}

To prepare for \cref{sec:strict_sbf,sec:applications}, this section reviews permutative categories \pcref{def:symmoncat}, symmetric monoidal functors \pcref{def:monoidalfunctor}, monoidal natural transformations \pcref{def:monnat}, bipermutative categories \pcref{def:bipermcat}, Laplaza's Coherence Theorem \pcref{thm:laplaza}, symmetric bimonoidal functors \pcref{def:sbfunctor}, and bimonoidal natural transformations \pcref{def:bimonnat}.  A detailed reference for this section is \cite[Part 1]{cerberusI}.  

\subsection*{Symmetric Monoidal Categories}

This subsection recalls the definitions of (symmetric) monoidal categories, permutative categories, (symmetric) monoidal functors, and monoidal natural transformations.  A reference for this material is \cite[Ch.\! 7 and 11]{maclane}.

\begin{definition}\label{def:monoidalcategory}
A \emph{monoidal category} $(\A,\otimes,\tu,\alpha,\lambda,\rho)$ consists of a category $\A$, a functor
\[\otimes \cn \A \times \A \to \A\]
called the \emph{monoidal product}, an object $\tu \in \A$ called the \emph{monoidal unit}, and natural isomorphisms for objects $a,b,c \in \A$ as follows.  
\begin{equation}\label{moncat_constraints}
\begin{tikzcd}[column sep=large]
(a \otimes b) \otimes c \ar{r}{\alpha_{a,b,c}}[swap]{\iso} & a \otimes (b \otimes c)
\end{tikzcd}
\qquad
\begin{tikzcd}[column sep=large]
\tu \otimes a \ar{r}{\lambda_a}[swap]{\iso} & a & a \otimes \tu \ar{l}{\iso}[swap]{\rho_a}
\end{tikzcd}
\end{equation}
They are called the \emph{associativity isomorphism}, the \emph{left unit isomorphism}, and the \emph{right unit isomorphism}.  These data are required to render the following unity and pentagon diagrams commutative for objects $a,b,c,d \in \A$.
\[\begin{tikzpicture}[xscale=1,yscale=1,vcenter]
\tikzset{0cell/.append style={nodes={scale=.8}}}
\tikzset{1cell/.append style={nodes={scale=.75}}}
\def\e{2.5} \def\h{1.8} \def\g{1.7} \def\v{.8} \def\u{2} \def\m{1em}
\draw[0cell]
(0,0) node (a) {(a \otimes \tu) \otimes b}
(a)++(\e,0) node (b) {a \otimes (\tu \otimes b)}
(a)++(\e/2,-1) node (c) {a \otimes b}
;
\draw[1cell]  
(a) edge node[swap,pos=.3] {\rho_a \otimes 1_b} (c)
(a) edge node {\alpha_{a,\tu,b}} (b)
(b) edge node[pos=.3] {1_a \otimes \lambda_b} (c)
;
\begin{scope}[shift={(7,.5)}]
\draw[0cell]
(0,0) node (x0) {(a\otimes b)\otimes (c\otimes d)}
(x0)++(-\h,-\v) node (x11) {((a\otimes b)\otimes c)\otimes d}
(x0)++(\h,-\v) node (x12) {a\otimes (b\otimes (c\otimes d))}
(x0)++(-\g,-\u) node (x21) {(a\otimes (b\otimes c)) \otimes d}
(x0)++(\g,-\u) node (x22) {a\otimes ((b\otimes c)\otimes d)}
;
\draw[1cell]
(x11) edge[transform canvas={xshift=-\m}] node[pos=.2] {\al_{a\otimes b,c,d}} (x0)
(x0) edge[transform canvas={xshift=\m}] node[pos=.8] {\al_{a,b,c\otimes d}} (x12)
(x11) edge node[swap,pos=.25] {\al_{a,b,c} \otimes 1_d} (x21)
(x21) edge node[inner sep=2pt] {\al_{a,b\otimes c,d}} (x22)
(x22) edge node[swap,pos=.7] {1_a \otimes \al_{b,c,d}} (x12)
;
\end{scope}
\end{tikzpicture}\]
A monoidal category is \emph{strict} if $\alpha$, $\lambda$, and $\rho$ are identities.
\end{definition}

\begin{definition}\label{def:symmoncat}
A \emph{symmetric monoidal category} $(\A,\brd)$ consists of a monoidal category $(\A,\otimes,\tu,\al,\la,\rho)$ and a natural isomorphism
\[a \otimes b \fto[\iso]{\xi_{a,b}} b \otimes a \forspace a,b \in \A,\]
called the \emph{braiding}, such that the following symmetry and hexagon diagrams commute for $a,b,c \in \A$. 
\begin{equation}\label{smc_axioms}
\begin{tikzpicture}[xscale=3,yscale=1,vcenter]
\def\h{.1}
\draw[0cell=.9] 
(0,0) node (a) {a \otimes b}
(a)++(.7,0) node (c) {a \otimes b}
(a)++(.35,-1) node (b) {b \otimes a}
;
\draw[1cell=.8] 
(a) edge node {1_{a \otimes b}} (c)
(a) edge node [swap,pos=.3] {\brd_{a,b}} (b)
(b) edge node [swap,pos=.7] {\brd_{b,a}} (c)
;
\begin{scope}[shift={(1.8,.5)}]
\draw[0cell=.8] 
(0,0) node (x11) {(b \otimes a) \otimes c}
(x11)++(.9,0) node (x12) {b \otimes (a \otimes c)}
(x11)++(-\h,-1) node (x21) {(a \otimes b) \otimes c}
(x12)++(\h,-1) node (x22) {b \otimes (c \otimes a)}
(x11)++(0,-2) node (x31) {a \otimes (b \otimes c)}
(x12)++(0,-2) node (x32) {(b \otimes c) \otimes a}
;
\draw[1cell=.8]
(x21) edge node[pos=.25] {\brd_{a,b} \otimes 1_c} (x11)
(x11) edge node {\alpha_{b,a,c}} (x12)
(x12) edge node[pos=.75] {1_b \otimes \brd_{a,c}} (x22)
(x21) edge node[swap,pos=.25] {\alpha_{a,b,c}} (x31)
(x31) edge node {\brd_{a,b \otimes c}} (x32)
(x32) edge node[swap,pos=.7] {\alpha_{b,c,a}} (x22)
;
\end{scope}
\end{tikzpicture}
\end{equation}
A \emph{permutative category} is a strict symmetric monoidal category, which means that $\al$, $\la$, and $\rho$ are identities.
\end{definition}

\begin{definition}\label{def:monoidalfunctor}
Suppose $\A$ and $\B$ are monoidal categories.  A \emph{monoidal functor} 
\[(\sff, \sff^2, \sff^0) \cn \A \to \B\]
consists of a functor $\sff \cn \A \to \B$, a morphism $\sff^0 \cn \tu \to \sff\tu$ called the \emph{unit constraint}, and a natural transformation
\begin{equation}\label{mon_constraint}
\begin{tikzcd}[column sep=large]
\sff a \otimes \sff b \ar{r}{\sff^2_{a,b}} & \sff(a \otimes b)
\end{tikzcd}
\end{equation}
for $a,b \in \A$, called the \emph{monoidal constraint}, such that the following unity and associativity diagrams commute for objects $a,b,c \in \A$.  
\begin{equation}\label{mfaxioms}
\begin{tikzcd}
\tu \otimes \sff a \dar[swap]{\sff^0 \otimes 1_{\sff a}} \rar{\lambda_{\sff a}} & \sff a \\ 
\sff\tu \otimes \sff a \rar{\sff^2_{\tu,a}} & \sff(\tu \otimes a) \uar[swap]{\sff\lambda_a} \\
\sff a \otimes \tu \dar[swap]{1_{\sff a} \otimes \sff^0} \rar{\rho_{\sff a}} & \sff a \\ 
\sff a \otimes \sff\tu \rar{\sff^2_{a,\tu}} & \sff(a \otimes \tu) \uar[swap]{\sff\rho_a}
\end{tikzcd}
\qquad
\begin{tikzcd}[column sep=large]
\bigl(\sff a \otimes \sff b\bigr) \otimes \sff c \rar{\alpha} \dar[swap]{\sff^2_{a,b} \otimes 1_{\sff c}} 
& \sff a \otimes \bigl(\sff b \otimes \sff c\bigr) \dar{1_{\sff a} \otimes \sff^2_{b,c}}\\
\sff(a \otimes b) \otimes \sff c \dar[swap]{\sff^2_{a \otimes b,c}} & \sff a \otimes \sff(b \otimes c) \dar[d]{\sff^2_{a,b \otimes c}}\\
\sff\bigl((a \otimes b) \otimes c\bigr) \rar{\sff\alpha} & \sff\bigl(a \otimes (b \otimes c)\bigr)
\end{tikzcd}
\end{equation}
A monoidal functor is \emph{strictly unital} if $\sff^0$ is the identity; \emph{strong} if $\sff^0$ and $\sff^2$ are isomorphisms; and \emph{strict} if $\sff^0$ and $\sff^2$ are identities.  

Suppose $(\sfh,\sfh^2,\sfh^0) \cn \B \to \C$ is another monoidal functor.  The \emph{composite} monoidal functor
\begin{equation}\label{monfunctor_comp}
\A \fto{(\sfh\sff, (\sfh\sff)^2, (\sfh\sff)^0)} \C
\end{equation}
is defined by the composite functor $\sfh\sff \cn \A \to \C$, the unit constraint
\[\tu \fto{\sfh^0} \sfh\tu \fto{\sfh(\sff^0)} \sfh\sff\tu,\]
and the monoidal constraint
\[(\sfh\sff a) \otimes (\sfh\sff b) 
\fto{\sfh^2_{\sff a, \sff b}} \sfh(\sff a \otimes \sff b) 
\fto{\sfh(\sff^2_{a,b})} \sfh\sff (a \otimes b)\]
for $a,b \in \A$.

A monoidal functor between symmetric monoidal categories is a \emph{symmetric monoidal functor} if the diagram
\begin{equation}\label{ftwobraiding}
\begin{tikzcd}[column sep=large]
\sff a \otimes \sff b \dar[swap]{\sff^2_{a,b}} \rar{\brd_{\sff a,\sff b}} & \sff b \otimes \sff a \dar{\sff^2_{b,a}} \\ 
\sff(a \otimes b) \rar{\sff\brd_{a,b}} & \sff(b \otimes a)
\end{tikzcd}
\end{equation}
commutes for $a,b \in \A$.  The category of small permutative categories and symmetric monoidal functors is denoted by $\perm$.  Its subcategory of strict symmetric monoidal functors is denoted by $\permst$.
\end{definition}

\begin{definition}\label{def:monnat}
Suppose $(\sff,\sff^2,\sff^0), (\sfh,\sfh^2,\sfh^0) \cn \A \to \B$ are monoidal functors between monoidal categories.  A \emph{monoidal natural transformation} $\uptha \cn \sff \to \sfh$ is a natural transformation between the underlying functors such that the following two diagrams in $\B$ commute for objects $a,b \in \A$.
\begin{equation}\label{monnat_axioms}
\begin{tikzpicture}[baseline={(a1.base)}]
\def\v{1.3}
\draw[0cell]
(0,0) node (a1) {\tu}
(a1)++(1.4,\v/2) node (a2) {\sff\tu}
(a2)++(0,-\v) node (a3) {\sfh\tu}
;
\draw[1cell=.9]
(a1) edge node[pos=.6] {\sff^0} (a2)
(a1) edge node[swap,pos=.6] {\sfh^0} (a3)
(a2) edge node {\uptha_\tu} (a3)
;
\begin{scope}[shift={(4,\v/2)}]
\draw[0cell]
(0,0) node (a11) {\sff a \otimes \sff b}
(a11)++(2.75,0) node (a12) {\sff(a \otimes b)}
(a11)++(0,-\v) node (a21) {\sfh a \otimes \sfh b}
(a12)++(0,-\v) node (a22) {\sfh(a \otimes b)}
;
\draw[1cell=.9]
(a11) edge node {\sff^2_{a,b}} (a12)
(a12) edge node {\uptha_{a \otimes b}} (a22)
(a11) edge node[swap] {\uptha_a \otimes \uptha_b} (a21)
(a21) edge node {\sfh^2_{a,b}} (a22)
;
\end{scope}
\end{tikzpicture}
\end{equation}
A \emph{monoidal adjunction} between monoidal categories is an adjunction such that both adjoints are monoidal functors, and the unit and the counit are monoidal natural transformations.
\end{definition}

\subsection*{Bipermutative Categories}

Bipermutative categories have two permutative structures, one additive and one multiplicative, that are related by a pair of distributivity morphisms.  The original definition is given in \cite[page 154]{may-einfinity}.  The definition in the following form is given in \cite[2.5.2]{cerberusI}.  

\begin{definition}\label{def:bipermcat}
A \emph{bipermutative category} is a tuple
\[\big(\A,(\bplus,\bzero,\abrd),(\btimes,\bunit,\mbrd),(\ladot,\rhodot),(\ldist,\rdist)\big)\]
consisting of the following data.
\begin{itemize}
\item $\Aplus = (\A,\bplus,\bzero,\abrd)$ is a permutative category \pcref{def:symmoncat}, called the \emph{additive structure}.
\item $\Atimes = (\A,\btimes,\bunit,\mbrd)$ is a permutative category, called the \emph{multiplicative structure}.
\item $\ladot$, $\rhodot$, $\mbrd_{-,\bzero}$, and $\rdist$ are identity natural transformations, and $\ldist$ is a natural isomorphism, for objects $a,b,c \in \A$ as follows.
\begin{equation}\label{larhordist}
\begin{tikzpicture}[xscale=1,yscale=1,vcenter]
\def\h{2} \def\g{3.5} \def\v{-.9} \def\a{1.8}
\draw[0cell=1] 
(0,0) node (a1) {\bzero \btimes a}
(a1)++(\h,0) node (a2) {\bzero}
(a2)++(\h,0) node (a3) {a \btimes \bzero}
(a3)++(\h,0) node (a4) {\bzero \btimes a}
(a1)++(-\a,\v) node (b1) {(a \bplus b) \btimes c}
(b1)++(\g,0) node (b2) {(a \btimes c) \oplus (b \btimes c)}
(a4)++(\a,\v) node (c2) {(a \btimes b) \bplus (a \btimes c)}
(c2)++(-\g,0) node (c1) {a \btimes (b \bplus c)}
;
\draw[1cell=.85] 
(a1) edge node {\ladot_a} node[swap]{=} (a2)
(a3) edge node [swap] {\rhodot_a}  node {=} (a2)
(a3) edge node {\mbrd_{a,\bzero}} node[swap]{=} (a4)
(b1) edge node {\rdist_{a,b,c}} node[swap]{=} (b2)
(c1) edge node {\ldist_{a,b,c}} node[swap]{\iso} (c2)
;
\end{tikzpicture}
\end{equation}
\end{itemize}
These data are required to make the following three diagrams commute for $a,b,c,d \in \A$, where $\btimes$ is denoted by concatenation.  In the absence of parentheses, $\btimes$ takes precedence over $\bplus$.
\begin{equation}\label{biperm-axioms}
\begin{tikzpicture}[xscale=1,yscale=1,vcenter]
\def\h{2.4} \def\v{-1.2} \def\g{3.7} \def\m{1.5ex}
\draw[0cell=.9]
(0,0) node (a11) {a(b \bplus c)}
(a11)++(\h,0) node (a12) {ab \bplus ac}
(a11)++(0,\v) node (a21) {(b \bplus c)a}
(a12)++(0,\v) node (a22) {ba \bplus ca}
;
\draw[1cell=.8]
(a11) edge node {\ldist} (a12)
(a11) edge node[swap] {\mbrd} (a21)
(a21) edge node {\rdist} node[swap]{=} (a22)
(a22) edge node[swap] {\mbrd \bplus \mbrd} (a12)
;
\begin{scope}[shift={(0,-2)}]
\draw[0cell=.9]
(0,0) node (a11) {(a \bplus b)c}
(a11)++(\h,0) node (a12) {ac \bplus bc}
(a11)++(0,\v) node (a21) {(b \bplus a)c}
(a12)++(0,\v) node (a22) {bc \bplus ac}
;
\draw[1cell=.8]
(a11) edge node {\rdist} node[swap]{=} (a12)
(a12) edge node {\abrd} (a22)
(a11) edge node[swap] {\abrd 1} (a21)
(a21) edge node {\rdist} node[swap]{=} (a22)
;
\end{scope}
\begin{scope}[shift={(5.5,-.5)}]
\draw[0cell=.9]
(0,0) node (x11) {(a \bplus b)(c \bplus d)}
(x11)++(\g,0) node (x12) {a(c \bplus d) \bplus b(c \bplus d)}
(x11)++(0,\v) node (x21) {(a \bplus b)c \bplus (a \bplus b)d} 
(x21)++(\g,0) node (x22) {ac \bplus ad \bplus bc \oplus bd}
(x21)++(\g/2,-1) node (x31) {ac \bplus bc \bplus ad \bplus bd}
;
\draw[1cell=.8]
(x11) edge node {\rdist} node[swap]{=} (x12)
(x12) edge node {\ldist \bplus \ldist} (x22)
(x22) edge[transform canvas={xshift=\m}] node[pos=.3] {1 \bplus \abrd \bplus 1} (x31)
(x11) edge node[swap] {\ldist} (x21)
(x21) edge[transform canvas={xshift=-\m}] node[sloped,pos=.6,inner sep=2pt]{=} node[swap,pos=.3] {\rdist \bplus \rdist} (x31)
; 
\end{scope}
\end{tikzpicture}
\end{equation}
The natural transformations $\ladot$, $\rhodot$, $\ldist$, and $\rdist$ are called, respectively, the \emph{left multiplicative zero}, the \emph{right multiplicative zero}, the \emph{left distributivity}, and the \emph{right distributivity}.  
\end{definition}

\begin{example}[Pointed Finite Sets]\label{ex:Fsk_dist}
There is a bipermutative category $(\Fsk,\wed,\sma)$ given by the category $\Fsk$ whose objects are pointed finite sets $\ordn = \{0,1,\ldots,n\}$ for $n \geq 0$ with basepoint 0.  Morphisms are pointed functions.
\begin{itemize}
\item The additive structure is given by the wedge product
\[\ordm \wed \ordn = \{0,1,\ldots,m,m+1,\ldots,m+n\}.\]
The additive monoidal unit is $\ord{0}$.  The additive braiding 
\[\ordm \wed \ordn \fto[\iso]{\abrd} \ordn \wed \ordm\]
is the block permutation that sends $\{1,\ldots,m\}$ to $\{n+1,\ldots,n+m\}$ and $\{m+1,\ldots,m+n\}$ to $\{1,\ldots,n\}$.
\item The multiplicative structure is given by the smash product
\[\ordm \sma \ordn = \ord{mn},\]
which uses the lexicographic ordering to identify a pair 
\begin{equation}\label{lex_order}
(i,j) \in \ordm \sma \ordn \withspace j+(i-1)n \in \ord{mn}
\end{equation} 
for $i,j > 0$.  If either $i=0$ or $j=0$, then $(i,j)$ is identified with the basepoint $0 \in \ord{mn}$.  The multiplicative monoidal unit is $\ord{1} = \{0,1\}$.  The multiplicative braiding
\[\ordm \sma \ordn \fto[\iso]{\mbrd} \ordn \sma \ordm\]
sends each non-basepoint $(i,j) \in \ordm \sma \ordn$ to $(j,i) \in \ordn \sma \ordm$.
\end{itemize}
The right distributivity $\rdist$ is the identity because $\ordm \sma \ordn = \ord{mn}$ is defined by the lexicographic ordering.  The left distributivity $\ldist$ is a non-identity isomorphism in general.
\end{example}

Bipermutative categories can be obtained by strictifying symmetric bimonoidal categories with invertible distributivity morphisms; see \cite[Prop.\! 6.3.5]{may-einfinity} for the original statement and \cite[5.4.7]{cerberusI} for a detailed proof.  Analogous to the situation for monoidal categories, this strictification theorem is a consequence of a coherence theorem for symmetric bimonoidal categories due to Laplaza \cite{laplaza,laplaza2}.  A detailed proof of this coherence theorem, along with corrections of Laplaza's original proof, is given in \cite[3.9.1 and 4.5.8]{cerberusI}.  The essence of Laplaza's Coherence Theorem is that, for a symmetric bimonoidal category whose distributivity morphisms satisfy a monomorphism condition, each well-defined formal diagram commutes.  Bipermutative categories are examples of symmetric bimonoidal categories with invertible distributivity morphisms, so Laplaza's Coherence Theorem also applies to bipermutative categories.  Since this paper only uses Laplaza's Coherence Theorem for bipermutative categories, it is stated in the following simplified form.  The reader is referred to \cite[Ch.\! 3 and 4]{cerberusI} for more detail.

\begin{theorem}[Laplaza]\label{thm:laplaza}
Each well-defined formal diagram in a bipermutative category commutes.
\end{theorem}

\subsection*{Bimonoidal Functors}
Symmetric bimonoidal functors between symmetric bimonoidal categories are defined in \cite[5.1.1]{cerberusI}.  The restriction to bipermutative categories is given as follows.

\begin{definition}\label{def:sbfunctor}
Suppose $\A$ and $\B$ are bipermutative categories.  A \emph{symmetric bimonoidal functor}
\[\big(\sff,\ftwoplus,\fzeroplus,\ftwotimes,\fzerotimes\big) \cn \A \to \B\]
consists of two symmetric monoidal functors between permutative categories \pcref{def:monoidalfunctor}
\[\begin{split}
\Aplus = (\A, \bplus, \bzero, \abrd) & \fto{\fplus = (\sff,\ftwoplus,\fzeroplus)} \Bplus = (\B, \bplus, \bzero, \abrd) \andspace\\
\Atimes = (\A, \btimes, \bunit, \mbrd) & \fto{\ftimes = (\sff,\ftwotimes,\fzerotimes) } \Btimes = (\B, \btimes, \bunit, \mbrd),
\end{split}\]
called the \emph{additive structure} and the \emph{multiplicative structure}, such that the following two diagrams in $\B$ commute for objects $a,b,c \in \A$.  The product $\btimes$ is denoted by concatenation.  In the absence of parentheses, $\btimes$ takes precedence over $\bplus$
\begin{equation}\label{sbf-axioms}
\begin{tikzpicture}[xscale=1,yscale=1,vcenter]
\def\h{2.3} \def\v{-1.2} \def\g{3.4}
\draw[0cell=.9]
(0,0) node (a11) {(\sff a) \bzero}
(a11)++(\h,0) node (a12) {(\sff a)(\sff\bzero)}
(a12)++(0,\v) node (a22) {\sff(a\bzero)}
(a11)++(0,2*\v) node (a31) {\bzero}
(a31)++(\h,0) node (a32) {\sff\bzero}
;
\draw[1cell=.85]
(a11) edge node {1 \fzeroplus} (a12)
(a12) edge node {\ftwotimes} (a22)
(a22) edge node {\sff\rhodot} node[swap,sloped] {=} (a32)
(a11) edge node[swap] {\rhodot} node[sloped] {=} (a31)
(a31) edge node {\fzeroplus} (a32)
;
\begin{scope}[shift={(5,0)}]
\draw[0cell=.9]
(0,0) node (b11) {(\sff a \bplus \sff b)(\sff c)}
(b11)++(\g,0) node (b12) {(\sff a)(\sff c) \bplus (\sff b)(\sff c)}
(b11)++(0,\v) node (b21) {\big(\sff(a \bplus b)\big) (\sff c)}
(b12)++(0,\v) node (b22) {\sff(ac) \bplus \sff(bc)}
(b21)++(0,\v) node (b31) {\sff\big((a \bplus b)c\big)}
(b22)++(0,\v) node (b32) {\sff(ac \oplus bc)}
;
\draw[1cell=.85]
(b11) edge node {\rdist} node[swap] {=} (b12)
(b12) edge node {\ftwotimes \bplus \ftwotimes} (b22)
(b22) edge node {\ftwoplus} (b32)
(b11) edge node[swap] {\ftwoplus 1} (b21)
(b21) edge node[swap] {\ftwotimes} (b31)
(b31) edge node {\sff\rdist} node[swap] {=} (b32)
;
\end{scope}
\end{tikzpicture}
\end{equation}
Composition of symmetric bimonoidal functors is defined by composing the additive structures and composing the multiplicative structures \cref{monfunctor_comp}.  An identity symmetric bimonoidal functor is defined by an identity functor $\sff$ and identity $\ftwoplus$, $\fzeroplus$, $\ftwotimes$, and $\fzerotimes$.

A symmetric bimonoidal functor is
\begin{itemize}
\item \emph{strictly unital} if the additive and multiplicative unit constraints
\[\bzero \fto{\fzeroplus} \sff\bzero \andspace \bunit \fto{\fzerotimes} \sff\bunit\] 
are identities;
\item \emph{strong} if both $\fplus$ and $\ftimes$ are strong, which means that $\ftwoplus$, $\fzeroplus$, $\ftwotimes$, and $\fzerotimes$ are isomorphisms; and
\item \emph{strict} if both $\fplus$ and $\ftimes$ are strict, which means that $\ftwoplus$, $\fzeroplus$, $\ftwotimes$, and $\fzerotimes$ are identities.
\end{itemize}
When \cref{def:monoidalfunctor} applies to $\fplus$ or $\ftimes$, the words \emph{additively} or \emph{multiplicatively} are used.  A symmetric bimonoidal functor is \emph{multiplicatively strong} if $\ftimes$ is strong, which means that the multiplicative monoidal and unit constraints
\begin{equation}\label{mult_strong}
\sff a \btimes \sff b \fto[\iso]{\ftwotimes} \sff(a \btimes b) \andspace 
\bunit \fto[\iso]{\fzerotimes} \sff\bunit
\end{equation}
are isomorphisms.  

There are categories as follows.
\begin{itemize}
\item $\biperm$: small bipermutative categories and symmetric bimonoidal functors.
\item $\bipermms$: small bipermutative categories and multiplicatively strong symmetric bimonoidal functors.
\item $\bipermst$: small bipermutative categories and strict symmetric bimonoidal functors.\defmark
\end{itemize}
\end{definition}

We emphasize that for a multiplicatively strong symmetric bimonoidal functor $\sff$, the additive monoidal and unit constraints
\[\sff a \bplus \sff b \fto{\ftwoplus} \sff(a \bplus b) \andspace \bzero \fto{\fzeroplus} \sff\bzero\] 
are not required to be isomorphisms.

\begin{remark}\label{rk:sbf}
Consider \cref{def:sbfunctor}.
\begin{enumerate}
\item The two diagrams in \cref{sbf-axioms} are stated in terms of the right multiplicative zero $\rhodot$ and the right distributivity $\rdist$.  Those two diagrams are equivalent to the following diagrams stated in terms of the left multiplicative zero $\ladot$ and the left distributivity $\ldist$ in \cref{larhordist}.
\begin{equation}\label{sbf-axioms-left}
\begin{tikzpicture}[xscale=1,yscale=1,vcenter]
\def\h{2.3} \def\v{-1.2} \def\g{3.4}
\draw[0cell=.9]
(0,0) node (a11) {\bzero(\sff a)}
(a11)++(\h,0) node (a12) {(\sff\bzero)(\sff a)}
(a12)++(0,\v) node (a22) {\sff(\bzero a)}
(a11)++(0,2*\v) node (a31) {\bzero}
(a31)++(\h,0) node (a32) {\sff\bzero}
;
\draw[1cell=.85]
(a11) edge node {\fzeroplus 1} (a12)
(a12) edge node {\ftwotimes} (a22)
(a22) edge node {\sff\ladot} node[swap,sloped] {=} (a32)
(a11) edge node[swap] {\ladot} node[sloped] {=} (a31)
(a31) edge node {\fzeroplus} (a32)
;
\begin{scope}[shift={(5,0)}]
\draw[0cell=.9]
(0,0) node (b11) {(\sff a)(\sff b \bplus \sff c)}
(b11)++(\g,0) node (b12) {(\sff a)(\sff b) \bplus (\sff a)(\sff c)}
(b11)++(0,\v) node (b21) {(\sff a)\big(\sff(b \bplus c)\big)}
(b12)++(0,\v) node (b22) {\sff(ab) \bplus \sff(ac)}
(b21)++(0,\v) node (b31) {\sff\big(a(b \bplus c)\big)}
(b22)++(0,\v) node (b32) {\sff(ab \oplus ac)}
;
\draw[1cell=.85]
(b11) edge node {\ldist} (b12)
(b12) edge node {\ftwotimes \bplus \ftwotimes} (b22)
(b22) edge node {\ftwoplus} (b32)
(b11) edge node[swap] {1\ftwoplus} (b21)
(b21) edge node[swap] {\ftwotimes} (b31)
(b31) edge node {\sff\ldist} (b32)
;
\end{scope}
\end{tikzpicture}
\end{equation}
See \cite[5.1.4]{cerberusI} for a proof of this statement. 
\item In \cite[page 155]{may-einfinity}, a \emph{morphism of symmetric bimonoidal categories} between bipermutative categories is the same as our strictly unital strong symmetric bimonoidal functor. 
A \emph{morphism of bipermutative categories} in \cite[page 155]{may-einfinity} is the same as our strict symmetric bimonoidal functor.\defmark
\end{enumerate}
\end{remark}

\begin{definition}\label{def:bimonnat}
Given symmetric bimonoidal functors \pcref{def:sbfunctor}
\[\big(\sff,\ftwoplus,\fzeroplus,\ftwotimes,\fzerotimes\big) \andspace 
\big(\sfh,\htwoplus,\hzeroplus,\htwotimes,\hzerotimes\big) \cn \A \to \B,\]
a \emph{bimonoidal natural transformation} $\uptha \cn \sff \to \sfh$ is a natural transformation between the underlying functors that is a monoidal natural transformation \pcref{def:monnat} with respect to both the additive structures and the multiplicative structures.  A \emph{bimonoidal adjunction} between bipermutative categories is an adjunction such that both adjoints are symmetric bimonoidal functors, and the unit and the counit are bimonoidal natural transformations.
\end{definition}

Unpacking \cref{def:bimonnat}, a bimonoidal natural transformation $\uptha \cn \sff \to \sfh$ is a natural transformation that makes the following four diagrams in $\B$ commute for objects $a,b \in \A$.
\begin{equation}\label{bimonnat_axioms}
\begin{tikzpicture}[vcenter]
\def\v{1.3} \def\s{4.5} \def\u{-2.3}
\draw[0cell]
(0,0) node (a1) {\bzero}
(a1)++(1.4,\v/2) node (a2) {\sff\bzero}
(a2)++(0,-\v) node (a3) {\sfh\bzero}
;
\draw[1cell=.9]
(a1) edge node[pos=.6] {\fzeroplus} (a2)
(a1) edge node[swap,pos=.6] {\hzeroplus} (a3)
(a2) edge node {\uptha_\bzero} (a3)
;
\begin{scope}[shift={(\s,\v/2)}]
\draw[0cell]
(0,0) node (a11) {\sff a \bplus \sff b}
(a11)++(2.75,0) node (a12) {\sff(a \bplus b)}
(a11)++(0,-\v) node (a21) {\sfh a \bplus \sfh b}
(a12)++(0,-\v) node (a22) {\sfh(a \bplus b)}
;
\draw[1cell=.9]
(a11) edge node {\ftwoplus} (a12)
(a12) edge node {\uptha_{a \bplus b}} (a22)
(a11) edge node[swap] {\uptha_a \bplus \uptha_b} (a21)
(a21) edge node {\htwoplus} (a22)
;
\end{scope}
\begin{scope}[shift={(0,\u)}]
\draw[0cell]
(0,0) node (a1) {\bunit}
(a1)++(1.4,\v/2) node (a2) {\sff\bunit}
(a2)++(0,-\v) node (a3) {\sfh\bunit}
;
\draw[1cell=.9]
(a1) edge node[pos=.6] {\fzerotimes} (a2)
(a1) edge node[swap,pos=.6] {\hzerotimes} (a3)
(a2) edge node {\uptha_\bunit} (a3)
;
\end{scope}
\begin{scope}[shift={(\s,\u+\v/2)}]
\draw[0cell]
(0,0) node (a11) {\sff a \btimes \sff b}
(a11)++(2.75,0) node (a12) {\sff(a \btimes b)}
(a11)++(0,-\v) node (a21) {\sfh a \btimes \sfh b}
(a12)++(0,-\v) node (a22) {\sfh(a \btimes b)}
;
\draw[1cell=.9]
(a11) edge node {\ftwotimes} (a12)
(a12) edge node {\uptha_{a \btimes b}} (a22)
(a11) edge node[swap] {\uptha_a \btimes \uptha_b} (a21)
(a21) edge node {\htwotimes} (a22)
;
\end{scope}
\end{tikzpicture}
\end{equation}

\section{Strictification of Multiplicatively Strong Symmetric Bimonoidal Functors}\label{sec:strict_sbf}

This section states and proves the main result of this paper, \cref{thm:sbfstrict}, which is a precise version of \cref{thm:main} in the \nameref{sec:introduction}.  \cref{rk:Bs_adjoint} discusses some subtle differences between \cref{may4.3,thm:sbfstrict}.

\begin{theorem}\label{thm:sbfstrict}
There is a functor
\[\bipermms \fto{\Bs} \bipermst\]
from the category $\bipermms$ of small bipermutative categories and multiplicatively strong symmetric bimonoidal functors to the category $\bipermst$ of small bipermutative categories and strict symmetric bimonoidal functors.  Moreover, for each bipermutative category $\A$, there is a bimonoidal adjunction
\[\begin{tikzpicture}
\def\b{22} \def\h{1.8}
\draw[0cell]
(0,0) node (a1) {\phantom{\A}}
(a1)++(-.15,0) node (a1') {\Bs\A}
(a1)++(\h,0) node (a2) {\A}
(a1)++(\h/2,0) node (a0) {\bot}
;
\draw[1cell=.9]
(a1) edge[bend left=\b] node {\lefta} (a2)
(a2) edge[bend left=\b] node {\righta} (a1)
;
\end{tikzpicture}\]
such that the following three statements hold.
\begin{enumerate}
\item\label{thm:sbfstrict_i} The counit $\counita \cn \lefta\righta \to 1_{\A}$ is the identity bimonoidal natural transformation.
\item\label{thm:sbfstrict_ii} The left adjoint $\lefta$ is an additively strict and multiplicatively strictly unital strong symmetric bimonoidal functor.  It is natural in strict symmetric bimonoidal functors.
\item\label{thm:sbfstrict_iii} The right adjoint $\righta$ is a multiplicatively strong symmetric bimonoidal functor.  It is natural in multiplicatively strong symmetric bimonoidal functors.  
\end{enumerate}
\end{theorem}

The application to multiplicative infinite loop space theory in \cref{sec:applications} uses only the statement of \cref{thm:sbfstrict} and not its proof.  The reader can safely jump to \cref{sec:applications} first and read the proof of \cref{thm:sbfstrict} later.

The rest of this section constructs the functors
\begin{itemize}
\item $\Bs$ in  \cref{def:BAobject,def:BAmorphisms,def:BAcomp,def:BA_mor_sumprod,def:BA_braidings,def:BA_mzero_dist,def:Bsf};
\item $\lefta$ in \cref{def:lefta}; and 
\item $\righta$ in \cref{def:righta}.
\end{itemize}
The proof of \cref{thm:sbfstrict} is given near the end of this section.  The naturality assertions about $\lefta$ and $\righta$ in \cref{thm:sbfstrict} \eqref{thm:sbfstrict_ii} and \eqref{thm:sbfstrict_iii} are explained in \cref{lefta_nat_diagram,righta_nat_diagram}.  The reason that the domain of $\Bs$ has \emph{multiplicatively strong} symmetric bimonoidal functors is the arrow $\ftwotimesinv$ in the diagram \cref{Bsf_mor_ii}, which is part of the definition of the strictified symmetric bimonoidal functor $\Bsf$ on morphisms.

\begin{assumption}\label{assu:sbf_strict}
For the rest of this section, $\A$ denotes an arbitrary bipermutative category \pcref{def:bipermcat}, with additive structure $\Aplus = (\A,\bplus,\bzero,\abrd)$ and multiplicative structure $\Atimes = (\A,\btimes,\bunit,\mbrd)$.  Smallness is only needed for the existence of the categories $\bipermms$ and $\bipermst$.
\end{assumption}

\subsection*{Object Assignment of $\Bs$}
To construct the functor $\Bs$ on objects, we need to define the bipermutative category $\BsA$.  To avoid confusion with the structures of $\A$, we denote the additive structure of $\BsA$ by $\BsAplus$ and the multiplicative structure by $\BsAtimes$.  We first define the object parts of $\BsAplus$ and $\BsAtimes$.  We denote the ordered set with $r \geq 0$ elements by 
\begin{equation}\label{ufs}
\ufs{r} = \{1 < 2 < \cdots < r\}
\end{equation}
with $\ufs{0} = \emptyset$.  The notation $i \in \ufs{r}$ means that the index $i$ runs through the finite set $\ufs{r}$, from 1 to $r$.  

\begin{definition}[Objects of $\BsA$]\label{def:BAobject}
We define objects of $\BsA$, the additive monoidal unit $\szero$, the multiplicative monoidal unit $\sunit$, the sum $\splus$ of objects, and the product $\stimes$ of objects as follows.
\begin{description}
\item[Objects] An object in $\BsA$ is a finite double sequence
\begin{equation}\label{BA_object}
\begin{split}
\bba &= \ang{\ang{a^i_j}_{j \in \ufs{m}_i}}_{i \in \ufs{r}}\\
&= \ang{\ang{a^1_j}_{j \in \ufs{m}_1} , \ldots , \ang{a^r_j}_{j \in \ufs{m}_r}}
\end{split}
\end{equation}
with each $a^i_j$ an object in $\A$.  We call $r \geq 0$ the \emph{additive length} of $\bba$.  For each $1 \leq i \leq r$, the $m_i$-tuple
\[a^i = \ang{a^i_j}_{j \in \ufs{m}_i} = \ang{a^i_1 , \ldots , a^i_{m_i}}\]
is called the \emph{$i$-th monomial} of $\bba$, with \emph{multiplicative length} $m_i \geq 0$ and \emph{$j$-th alphabet} $a^i_j$.  To give some intuition, we think of an object $\bba \in \BsA$ as a formal polynomial 
\[\sum_{i=1}^r a^i_1 a^i_2 \cdots a^i_{m_i}.\]
The terminology above and the definitions below are motivated by this interpretation of $\bba$.
\item[Additive monoidal unit] We define the unique object 
\begin{equation}\label{szero}
\szero = \emptyset \in \BsA
\end{equation}
with additive length 0 and no monomials.
\item[Multiplicative monoidal unit] We define the object 
\begin{equation}\label{sunit}
\sunit = \ang{\emptyset} \in \BsA
\end{equation}
with additive length 1, whose only monomial $\emptyset$ has multiplicative length 0.
\item[Sum] Suppose $\bba \in \BsA$ is an object as defined in \cref{BA_object} and
\begin{equation}\label{bbb}
\bbb = \ang{\ang{b^k_l}_{l \in \ufsn_k}}_{k \in \ufss} \in \BsA
\end{equation}
is an object with additive length $s$.  For each $1 \leq k \leq s$, its $k$-th monomial
\[b^k = \ang{b^k_l}_{l \in \ufsn_k}\]
has multiplicative length $n_k$ and $l$-th alphabet $b^k_l \in \A$.   We define their \emph{sum} as the object
\begin{equation}\label{splus}
\bba \splus \bbb = \ang{\ang{a^i}_{i \in \ufsr} , \ang{b^k}_{k \in \ufss}} \in \BsA
\end{equation}
with additive length $r+s$, obtained by concatenating $\bba$ and $\bbb$.
\begin{itemize}
\item For $1 \leq i \leq r$, its $i$-th monomial is $a^i$. 
\item For $1 \leq k \leq s$, its $(r+k)$-th monomial is $b^k$.
\end{itemize}
Note that the sum $\splus$ is associative on objects, and the additive monoidal unit $\szero$ defined in \cref{szero} is a two-sided unit for the sum:
\[\bba \splus \szero = \bba = \szero \splus \bba.\]
\item[Product] For objects $\bba, \bbb \in \BsA$ as defined in \cref{BA_object,bbb}, we define their \emph{product} as the object
\begin{equation}\label{stimes}
\begin{split}
\bba \stimes \bbb &= \ang{\ang{a^i b^k}_{k \in \ufss}}_{i \in \ufsr} \in \BsA\\
&= \ang{\ang{a^1 b^k}_{k \in \ufss} , \ldots , \ang{a^r b^k}_{k \in \ufss}}
\end{split}
\end{equation}
with additive length $rs$.  For $1 \leq i \leq r$ and $1 \leq k \leq s$, its $(k + (i-1)s)$-th monomial is the concatenation of the $i$-th monomial $a^i$ of $\bba$ and the $k$-th monomial $b^k$ of $\bbb$, denoted
\begin{equation}\label{aibk}
a^i b^k = \ang{\ang{a^i_j}_{j \in \ufsm_i} , \ang{b^k_l}_{l \in \ufsn_k}},
\end{equation}
with multiplicative length $m_i + n_k$.
\begin{itemize}
\item For $1 \leq j \leq m_i$, its $j$-th alphabet is $a^i_j$.
\item For $1 \leq l \leq n_k$, its $(m_i + l)$-th alphabet is $b^k_l$.
\end{itemize}
Note that the product $\bba \stimes \bbb$ is obtained by concatenating the monomials of $\bba$ and $\bbb$ and ordering them lexicographically.  The product $\stimes$ is associative on objects, and the multiplicative monoidal unit $\sunit = \ang{\emptyset}$ defined in \cref{sunit} is a two-sided unit for the product:
\[\bba \stimes \sunit = \bba = \sunit \stimes \bba.\]
\end{description}
This finishes the definition.
\end{definition}

Next, we define morphisms of $\BsA$.  Recall from \cref{ufs} that $\ufsr$ denotes the set $\{1,2,\ldots,r\}$ with its natural ordering.

\begin{definition}[Morphisms of $\BsA$]\label{def:BAmorphisms}
We consider objects 
\[\bba = \ang{\ang{a^i_j}_{j \in \ufsm_i}}_{i \in \ufs{r}} \andspace 
\bbb = \ang{\ang{b^k_l}_{l \in \ufsn_k}}_{k \in \ufss} \in \BsA\] 
as defined in \cref{BA_object,bbb}, with additive lengths $r$ and $s$.  The $i$-th monomial $a^i = \ang{a^i_j}_{j \in \ufsm_i}$ of $\bba$ has multiplicative length $m_i$, and the $k$-th monomial $b^k = \ang{b^k_l}_{l \in \ufsn_k}$ of $\bbb$ has multiplicative length $n_k$.  A \emph{morphism} 
\begin{equation}\label{BA_morphism}
\bba \fto{(\uphi; g)} \bbb \inspace \BsA
\end{equation}
consists of the following data.
\begin{description}
\item[Reindexing function] $\uphi \cn \ufsr \to \ufss$ is a function.
\item[Component morphisms] $g = \ang{g^k}_{k \in \ufss}$ consists of, for each index $k \in \ufss$, a morphism
\begin{equation}\label{BA_mor_ii}
\bigbplus_{i \in \uphiinv(k)} \bigbtimes_{j \in \ufsm_i} a^i_j 
\fto{g^k} \bigbtimes_{l \in \ufsn_k} b^k_l \inspace \A.
\end{equation}
In \cref{BA_mor_ii}, the sum $\bplus$ is indexed by the ordered subset $\uphiinv(k) \subseteq \ufsr$.  By convention, an empty sum is the additive monoidal unit $\bzero \in \A$, and an empty product $\btimes$ is the multiplicative monoidal unit $\bunit \in \A$ \pcref{def:bipermcat}.
\end{description}
The \emph{identity morphism} $1_\bba \cn \bba \to \bba$ of an object $\bba \in \BsA$ is defined by the identity function $1_{\ufsr}$ and identity morphisms of the objects $\bigbtimes_{j \in \ufsm_i} a^i_j$ for $i \in \ufsr$.
\end{definition}

Next, we define composition of morphisms of $\BsA$.

\begin{definition}[Composition of $\BsA$]\label{def:BAcomp}
We consider two morphisms \cref{BA_morphism}
\[\bba \fto{(\uphi; g)} \bbb \fto{(\uppsi; h)} \bbc \inspace \BsA,\]
where $\bbc$ is an object
\begin{equation}\label{bbc}
\bbc = \ang{\ang{c^p_q}_{q \in \ufst_p}}_{p \in \ufsu} \in \BsA
\end{equation}
with additive length $u$ \cref{BA_object}.  For each $1 \leq p \leq u$, its $p$-th monomial $c^p = \ang{c^p_q}_{q \in \ufst_p}$ has multiplicative length $t_p$.  The \emph{composite morphism}
\begin{equation}\label{BA_mor_comp}
\bba \fto{(\uppsi; h) \circ (\uphi; g) = (\uppsi\uphi ; hg)} \bbc 
\end{equation}
is defined as follows.
\begin{description}
\item[Reindexing function] $\uppsi\uphi$ is the composite function
\begin{equation}\label{BA_comp_re}
\ufsr \fto{\uphi} \ufss \fto{\uppsi} \ufsu.
\end{equation}
\item[Component morphisms]
The component morphism of $hg = \ang{(hg)^p}_{p \in \ufsu}$ at an index $p \in \ufsu$ is defined as the following composite in $\A$.
\begin{equation}\label{BA_comp_hg}
\begin{tikzpicture}[vcenter]
\def\v{-1.6}
\draw[0cell=.85]
(0,0) node (a11) {\bigbplus_{i \in (\uppsi\uphi)^\inv(p)} \bigbtimes_{j \in \ufsm_i} a^i_j}
(a11)++(4.8,0) node (a12) {\bigbtimes_{q \in \ufst_p} c^p_q}
(a11)++(0,\v) node (a21) {\bigbplus_{k \in \uppsiinv(p)} \bigbplus_{i \in \uphiinv(k)} \bigbtimes_{j \in \ufsm_i} a^i_j}
(a12)++(0,\v) node (a22) {\bigbplus_{k \in \uppsiinv(p)} \bigbtimes_{l \in \ufsn_k} b^k_l}
;
\draw[1cell=.8]
(a11) edge node {(hg)^p} (a12)
(a11) edge node[swap] {\si_{\uphi, \uppsi; p}} node {\iso} (a21)
(a21) edge node {\bigbplus_{k \in \uppsiinv(p)} g^k} (a22)
(a22) edge node[swap] {h^p} (a12)
;
\end{tikzpicture}
\end{equation}
In \cref{BA_comp_hg}, $g^k$ and $h^p$ are component morphisms of $g$ and $h$, as defined in \cref{BA_mor_ii}.  The left vertical arrow $\si_{\uphi, \uppsi; p}$ is the unique coherence isomorphism \cite[Theorem XI.1]{maclane} in the additive structure $\Aplus = (\A, \bplus, \bzero, \abrd)$, which is a permutative category.  It permutes the objects $\bigbtimes_{j \in \ufsm_i} a^i_j$ and inserts a copy of the additive monoidal unit $\bzero \in \A$ for each index $k \in \uppsiinv(p)$ with $\uphiinv(k) = \emptyset$.  
\end{description}
Note that the coherence isomorphism $\si_{\uphi, \uppsi; p}$ in \cref{BA_comp_hg} is the identity morphism if the restriction of $\uphi$ to the ordered subset $(\uppsi\uphi)^\inv(p) \subseteq \ufsr$ is order-preserving.  For example, this is the case if $|\uppsiinv(p)| = 1$.  Thus, identity morphisms are two-sided units for composition.  Associativity of composition follows from the uniqueness and naturality of coherence isomorphisms in the permutative category $\Aplus$.  This finishes the definition of the category $\BsA$.
\end{definition}

Sums and products of objects in $\BsA$ are already defined in \cref{splus,stimes}.  Next, we extend those constructions to morphisms \pcref{def:BAmorphisms}.

\begin{definition}[Sums and Products of Morphisms]\label{def:BA_mor_sumprod}
We consider two morphisms 
\[\begin{split}
\ang{\ang{a^i_j}_{j \in \ufsm_i}}_{i \in \ufsr} = \bba & 
\fto{(\uphi; g)} \bbb = \ang{\ang{b^k_l}_{l \in \ufsn_k}}_{k \in \ufss}\andspace\\ 
\ang{\ang{c^p_q}_{q \in \ufst_p}}_{p \in \ufsu} = \bbc & 
\fto{(\uppsi; h)} \bbd = \ang{\ang{d^y_z}_{z \in \ufsv_y}}_{y \in \ufsw}
\end{split}\]
in $\BsA$ as defined in \cref{BA_morphism}.  We define their sum and product as follows.
\begin{description}
\item[Sum] 
We define the morphism
\begin{equation}\label{BA_mor_sum}
\ang{\ang{a^i}_{i \in \ufsr} , \ang{c^p}_{p \in \ufsu}} = \bba \splus \bbc 
\fto[=  (\uphi \sqcup \uppsi; (g,h))]{(\uphi; g) \splus (\uppsi; h)} 
\bbb \splus \bbd = \ang{\ang{b^k}_{k \in \ufss} , \ang{d^y}_{y \in \ufsw}} 
\end{equation}
in $\BsA$ whose reindexing function is the coproduct
\begin{equation}\label{BA_mor_sum_re}
\ufsr \sqcup \ufsu \fto{\uphi \sqcup \uppsi} \ufss \sqcup \ufsw.
\end{equation}
The sum $(\uphi; g) \splus (\uppsi; h)$ has component morphism
\begin{itemize}
\item $g^k$ for each index $k \in \ufss \subseteq \ufss \sqcup \ufsw$ and
\item $h^y$ for each index $y \in \ufsw \subseteq \ufss \sqcup \ufsw$.
\end{itemize}
The functoriality, associativity (for both objects and morphisms), and unity (with monoidal unit $\szero$ defined in \cref{szero}) of $\splus$ follow from the universal properties of coproducts.
\item[Product]  
We define the morphism
\begin{equation}\label{BA_mor_prod}
\ang{\ang{a^i c^p}_{p \in \ufsu}}_{i \in \ufsr} = \bba \stimes \bbc 
\fto[= \big(\uphi \times \uppsi; (g \btimes h)\deinv\big)]{(\uphi; g) \stimes (\uppsi; h)} 
\bbb \stimes \bbd = \ang{\ang{b^k d^y}_{y \in \ufsw}}_{k \in \ufss}
\end{equation}
in $\BsA$ whose reindexing function is the product
\begin{equation}\label{BA_mor_prod_re}
\ufs{ru} = \ufsr \times \ufsu \fto{\uphi \times \uppsi} \ufss \times \ufsw = \ufs{sw}
\end{equation}
with respect to the lexicographic ordering \cref{lex_order} in the domain and the codomain.  Recall from \cref{aibk} that $b^k d^y$ is the concatenated monomial with multiplicative length $n_k + v_y$.  The component morphism of $(\uphi; g) \stimes (\uppsi; h)$ at a pair of indices $(k,y) \in \ufss \times \ufsw$ is defined as the following composite in $\A$.
\begin{equation}\label{BA_mor_prod_ii}
\begin{tikzpicture}[vcenter]
\draw[0cell=.85]
(0,0) node (a1) {\bigbplus_{i \in \uphiinv(k)} \bigbplus_{p \in \uppsiinv(y)} \Big(\bigbtimes_{j \in \ufsm_i} a^i_j\Big) \otimes \Big(\bigbtimes_{q \in \ufst_p} c^p_q \Big)}
(a1)++(2.1,-1.1) node (a2) {\phantom{\Big(\bigbplus_{i \in \uphiinv(k)} \bigbtimes_{j \in \ufsm_i} a^i_j\Big) \btimes \Big(\bigbplus_{p \in \uppsiinv(y)} \bigbtimes_{q \in \ufst_p} c^p_q \Big)}}
(a2)++(0,-.1) node (a2') {\Big(\bigbplus_{i \in \uphiinv(k)} \bigbtimes_{j \in \ufsm_i} a^i_j\Big) \btimes \Big(\bigbplus_{p \in \uppsiinv(y)} \bigbtimes_{q \in \ufst_p} c^p_q \Big)}
(a1)++(5,0) node (a3) {\Big(\bigbtimes_{l \in \ufsn_k} b^k_l \Big) \btimes \Big(\bigbtimes_{z \in \ufsv_y} d^y_z \Big)}
(a1)++(-1,0) node (a1') {\phantom{\bigbtimes_{q \in \ufst_p}}}
(a3)++(0,0) node (a3') {\phantom{\bigbtimes_{q}}}
;
\draw[1cell=.8]
(a1') [rounded corners=2pt] |- node[swap,pos=.25] {\deinv_{\uphi,\uppsi; k,y}} node[pos=.25] {\iso} (a2)
;
\draw[1cell=.8]
(a2) [rounded corners=2pt] -| node[swap,pos=.75] {g^k \btimes h^y} (a3')
;
\end{tikzpicture}
\end{equation}
The first arrow $\deinv_{\uphi,\uppsi; k,y}$ is the unique Laplaza coherence isomorphism \pcref{thm:laplaza} in the bipermutative category $\A$.
\begin{itemize}
\item If both $\uphiinv(k)$ and $\uppsiinv(y)$ are nonempty, then $\deinv_{\uphi,\uppsi; k,y}$ consists of inverses of the right and left distributivity isomorphisms \cref{larhordist}, $\rdist$ and $\ldist$.  Note that, if $|\uppsiinv(y)| = 1$, then $\deinv_{\uphi,\uppsi; k,y}$ is the identity morphism, since $\rdist$ is the identity in the bipermutative category $\A$.
\item If $\uphiinv(k) = \emptyset$, then $\deinv_{\uphi,\uppsi; k,y}$ is the inverse of the left multiplicative zero $\ladot$ of $\A$.
\item If $\uppsiinv(y) = \emptyset$, then $\deinv_{\uphi,\uppsi; k,y}$ is the inverse of the right multiplicative zero $\rhodot$ of $\A$. 
\end{itemize}
In the second arrow, $g^k$ and $h^y$ are the indicated component morphisms of $g$ and $h$, as defined in \cref{BA_mor_ii}.  

The functoriality, associativity (for both objects and morphisms), and unity (with monoidal unit $\sunit$ defined in \cref{sunit}) of $\stimes$ follow from the uniqueness and naturality of Laplaza coherence isomorphisms \pcref{thm:laplaza}.
\end{description}
For each of $\splus$ and $\stimes$, defining the associativity and unit isomorphisms \cref{moncat_constraints} to be identities, we have constructed the strict monoidal categories $(\BsA,\splus,\szero)$ and $(\BsA,\stimes,\sunit)$.
\end{definition}

To define the additive structure $\BsAplus$ and the multiplicative structure $\BsAtimes$, next we define their braidings.

\begin{definition}[Braidings]\label{def:BA_braidings}
We consider objects
\[\bba = \ang{\ang{a^i_j}_{j \in \ufsm_i}}_{i \in \ufsr} \andspace
\bbb = \ang{\ang{b^k_l}_{l \in \ufsn_k}}_{k \in \ufss}\]
in $\BsA$ as defined in \cref{BA_object,bbb}.
\begin{description}
\item[Additive braiding]
We define the isomorphism 
\begin{equation}\label{BA_sabrd}
\ang{\ang{a^i}_{i \in \ufsr}, \ang{b^k}_{k \in \ufss}} = \bba \splus \bbb 
\fto[\iso]{\sabrd = (\uptau; \ang{1})} 
\bbb \splus \bba = \ang{\ang{b^k}_{k \in \ufss}, \ang{a^i}_{i \in \ufsr}}
\end{equation}
whose reindexing function is the block permutation 
\begin{equation}\label{BA_abrd_reindex}
\ufsr \sqcup \ufss\fto[\iso]{\uptau} \ufss \sqcup \ufsr
\end{equation}
that interchanges $\ufsr$ and $\ufss$.  The component morphism of $\sabrd$ is the identity morphism of the object
\begin{itemize} 
\item $\bigbtimes_{l \in \ufsn_k} b^k_l$ for each index $k \in \ufss \subseteq \ufss \sqcup \ufsr$ and
\item $\bigbtimes_{j \in \ufsm_i} a^i_j$ for each index $i \in \ufsr \subseteq \ufss \sqcup \ufsr$.
\end{itemize}
The naturality of $\sabrd$ in $\bba$ and $\bbb$ follows from \cref{def:BAcomp} and the naturality of the block permutation $\uptau$ with respect to functions on $\ufsr$ and $\ufss$.  The symmetry and hexagon diagrams in \cref{smc_axioms} commute for $\sabrd$ because the block permutations $\uptau$ have these properties and that each component morphism of $\sabrd$ is an identity morphism.
\item[Multiplicative braiding]
We define the isomorphism 
\begin{equation}\label{BA_smbrd}
\ang{\ang{a^i b^k}_{k \in \ufss}}_{i \in \ufsr} = \bba \stimes \bbb 
\fto[\iso]{\smbrd = (\tau; \ang{\mbrd})} 
\bbb \stimes \bba = \ang{\ang{b^k a^i}_{i \in \ufsr}}_{k \in \ufss}
\end{equation}
whose reindexing function is the swapping permutation
\begin{equation}\label{BA_mbrd_reindex}
\ufs{rs} = \ufsr \times \ufss \fto[\iso]{\tau} \ufss \times \ufsr = \ufs{sr}
\end{equation}
that sends each pair of indices $(i,k) \in \ufsr \times \ufss$ to $(k,i) \in \ufss \times \ufsr$.  Recall from \cref{aibk} that $a^i b^k$ is the concatenated monomial with multiplicative length $m_i + n_k$.  The component morphism of $\smbrd$ at a pair of indices $(k,i) \in \ufss \times \ufsr$ is the multiplicative braiding in $\A$
\begin{equation}\label{BA_mbrd_component}
\Big(\bigbtimes_{j \in \ufsm_i} a^i_j\Big) \btimes \Big(\bigbtimes_{l \in \ufsn_k} b^k_l\Big) 
\fto[\iso]{\mbrd}
\Big(\bigbtimes_{l \in \ufsn_k} b^k_l\Big) \btimes \Big(\bigbtimes_{j \in \ufsm_i} a^i_j\Big)
\end{equation}
that permutes the objects $\bigbtimes_{j \in \ufsm_i} a^i_j$ and $\bigbtimes_{l \in \ufsn_k} b^k_l$

The naturality of $\smbrd$ is proved in \cref{smbrd_natural} below.  The symmetry and hexagon diagrams in \cref{smc_axioms} commute for $\smbrd$ because the swapping permutations $\tau$ and the multiplicative braiding $\mbrd$ in $\A$ both have these properties.
\end{description}
Moreover, the multiplicative braiding
\[\bba \stimes \szero \fto{\smbrd} \szero \stimes \bba\]
is the identity morphism of the additive monoidal unit $\szero$ \cref{szero}, since both $\bba \stimes \szero$ and $\szero \stimes \bba$ have additive length 0.
\end{definition}

\begin{lemma}\label{smbrd_natural}
$\smbrd$ in \cref{BA_smbrd} is a natural transformation.
\end{lemma}

\begin{proof}
Using the notation in \cref{def:BA_mor_sumprod}, we need to show that the diagram
\begin{equation}\label{smbrd_nat_i}
\begin{tikzpicture}[vcenter]
\def\v{-1.3}
\draw[0cell]
(0,0) node (a11) {\bba \stimes \bbc}
(a11)++(2.3,0) node (a12) {\bbc \stimes \bba}
(a11)++(0,\v) node (a21) {\bbb \stimes \bbd}
(a12)++(0,\v) node (a22) {\bbd \stimes \bbb}
;
\draw[1cell=.9]
(a11) edge node {\smbrd} (a12)
(a12) edge node {(\uppsi;h) \stimes (\uphi;g)} (a22)
(a11) edge node[swap] {(\uphi;g) \stimes (\uppsi;h)} (a21)
(a21) edge node {\smbrd} (a22)
;
\end{tikzpicture}
\end{equation}
in $\BsA$ commutes.  By \cref{BA_comp_re,BA_mor_prod_re,BA_mbrd_reindex}, the reindexing functions of the two composites in the diagram \cref{smbrd_nat_i} are given by the two composites in the following commutative diagram of functions.
\[\begin{tikzpicture}
\def\v{-1.3}
\draw[0cell]
(0,0) node (a11) {\ufsr \times \ufsu}
(a11)++(2.3,0) node (a12) {\ufsu \times \ufsr}
(a11)++(0,\v) node (a21) {\ufss \times \ufsw}
(a12)++(0,\v) node (a22) {\ufsw \times \ufss}
;
\draw[1cell=.9]
(a11) edge node {\tau} (a12)
(a12) edge node {\uppsi \times \uphi} (a22)
(a11) edge node[swap] {\uphi \times \uppsi} (a21)
(a21) edge node {\tau} (a22)
;
\end{tikzpicture}\]

By \cref{BA_comp_hg,BA_mor_prod_ii,BA_mbrd_component}, for each pair of indices $(y,k) \in \ufsw \times \ufss$, the $(y,k)$-component morphisms of the two composites in the diagram \cref{smbrd_nat_i} are given by the two boundary composites in the diagram \cref{smbrd_nat_ii} in $\A$ below, where we use the following abbreviations.
\[\begin{aligned}
a^i_\crdot &= \dbigbtimes_{j \in \ufsm_i} a^i_j & 
b^k_\crdot &= \dbigbtimes_{l \in \ufsn_k} b^k_l &
c^p_\crdot &= \dbigbtimes_{q \in \ufst_p} c^p_q \phantom{M} &
d^y_\crdot &= \dbigbtimes_{z \in \ufsv_y} d^y_z\\
\dbigbplus_i &= \dbigbplus_{i \in \uphiinv(k)} & 
\dbigbplus_p &= \dbigbplus_{p \in \uppsiinv(y)}
\end{aligned}\]
\begin{equation}\label{smbrd_nat_ii}
\begin{tikzpicture}[vcenter]
\def\v{1.3} \def\h{4} \def\g{1.5} \def\u{-1.1}
\draw[0cell=.9]
(0,0) node (a11) {\dbigbplus_i \, \dbigbplus_p (a^i_\crdot \btimes c^p_\crdot)}
(a11)++(0,\v) node (a12) {\dbigbplus_p \, \dbigbplus_i (a^i_\crdot \btimes c^p_\crdot)}
(a12)++(\h,0) node (a13) {\dbigbplus_p \, \dbigbplus_i (c^p_\crdot \btimes a^i_\crdot)}
(a13)++(\h,0) node (a14) {\big(\dbigbplus_p c^p_\crdot\big) \btimes \big(\dbigbplus_i a^i_\crdot\big)}
(a14)++(0,-\v) node (a15) {d^y_\crdot \btimes b^k_\crdot}
(a11)++(\g,\u) node (b1) {\big(\dbigbplus_i a^i_\crdot\big) \btimes \big(\dbigbplus_p c^p_\crdot\big)}
(a15)++(-\g,\u) node (b2) {b^k_\crdot \btimes d^y_\crdot}
;
\draw[1cell=.8]
(a11) edge node {\si} (a12)
(a12) edge node {\bplus_p \bplus_i \mbrd} (a13)
(a13) edge node {\deinv} (a14)
(a14) edge node {h^y \btimes g^k} (a15)
(a11) edge[transform canvas={xshift=-1ex}, shorten <=-1ex] node[swap,pos=.1] {\deinv} (b1)
(b1) edge node {g^k \btimes h^y} (b2)
(b2) edge[transform canvas={xshift=1ex}] node[swap,pos=.7] {\mbrd} (a15)
(b1) edge node {\mbrd} (a14)
;
\end{tikzpicture}
\end{equation}
In the diagram \cref{smbrd_nat_ii}, the upper left and lower right regions commute by, respectively, Laplaza's Coherence \cref{thm:laplaza} and the naturality of the multiplicative braiding $\mbrd$ in $\A$.
\end{proof}

This finishes the construction of the permutative categories 
\[\BsAplus = (\BsA, \splus, \szero, \sabrd) \andspace \BsAtimes = (\BsA, \stimes, \sunit, \smbrd),\]
which are the additive and multiplicative structures of $\BsA$.  To finish the construction of the bipermutative category $\BsA$, next we define its multiplicative zeros and distributivity morphisms \cref{larhordist}.

\begin{definition}[Multiplicative Zeros and Distributivity]\label{def:BA_mzero_dist}
Using the permutative categories $\BsAplus$ and $\BsAtimes$, we define identity natural transformations
\begin{equation}\label{BA_mzero_dist}
\begin{tikzpicture}[vcenter]
\def\h{2}
\draw[0cell]
(0,0) node (a11) {\szero \stimes \bba}
(a11)++(\h,0) node (a12) {\szero}
(a12)++(\h,0) node (a13) {\bba \stimes \szero}
(a11)++(0,-.8) node (b1) {(\bba \splus \bbb) \stimes \bbc}
(b1)++(3.5,0) node (b2) {(\bba \stimes \bbc) \splus (\bbb \stimes \bbc)}
;
\draw[1cell=.9]
(a11) edge node {\sladot} node[swap] {=} (a12)
(a13) edge node {=} node[swap] {\srhodot} (a12)
(b1) edge node {\srdist} node[swap] {=} (b2)
;
\end{tikzpicture}
\end{equation}
for objects $\bba, \bbb, \bbc \in \BsA$.
\begin{itemize}
\item The left multiplicative zero $\sladot$ and the right multiplicative zero $\srhodot$ are well defined and natural because $\szero \stimes \bba$ and $\bba \stimes \szero$ have additive length 0. 
\item The right distributivity $\srdist$ is well defined because, with the notation in \cref{BA_object,bbb,bbc}, both the domain and the codomain of $\srdist$ are equal to the object
\[\ang{\ang{\ang{a^i c^p}_{p \in \ufsu}}_{i \in \ufsr}, \ang{\ang{b^k c^p}_{p \in \ufsu}}_{k \in \ufss}} \in \BsA.\]
The naturality of $\srdist$ follows from \cref{def:BA_mor_sumprod}.
\end{itemize}

\parhead{Left distributivity}.  We define the left distributivity $\sldist$ for $\BsA$ using the upper left diagram in \cref{biperm-axioms}:
\begin{equation}\label{BA_ldist}
\begin{tikzpicture}[vcenter]
\def\v{-1.3} \def\u{.9} \def\h{1}
\draw[0cell]
(0,0) node (a11) {\bba \stimes (\bbb \splus \bbc)}
(a11)++(3.5,0) node (a12) {(\bba \stimes \bbb) \splus (\bba \stimes \bbc)}
(a11)++(0,\v) node (a21) {(\bbb \splus \bbc) \stimes \bba}
(a12)++(0,\v) node (a22) {(\bbb \stimes \bba) \splus (\bbc \stimes \bba)}
(a11)++(0,\u) node (a01) {\phantom{\ang{}_{p}}}
(a12)++(0,\u) node (a02) {\phantom{\ang{}_{p}}}
(a11)++(-\h,\u) node (a01') {\ang{\ang{a^i b^k}_{k \in \ufss}, \ang{a^i c^p}_{p \in \ufsu}}_{i \in \ufsr}}
(a12)++(\h,\u) node (a02') {\ang{\ang{\ang{a^i b^k}_{k \in \ufss}}_{i \in \ufsr} , \ang{\ang{a^i c^p}_{p \in \ufsu}}_{i \in \ufsr}}}
;
\draw[1cell=.9]
(a11) edge node {\sldist} node[swap] {\iso} (a12)
(a11) edge node[swap] {\smbrd} (a21)
(a21) edge node {\srdist} node[swap] {=} (a22)
(a22) edge node[swap] {\smbrd \splus \smbrd} (a12)
(a01) edge[equal] (a11)
(a02) edge[equal] (a12)
;
\end{tikzpicture}
\end{equation}
The naturality of $\sldist$ follows from the naturality of $\smbrd$ and $\srdist$, together with the functoriality of $\splus$ and $\stimes$ \pcref{def:BA_mor_sumprod,def:BA_braidings}.  

Using \cref{BA_mbrd_reindex} to unpack the diagram \cref{BA_ldist}, the reindexing function of $\sldist$ is the bijection
\begin{equation}\label{BA_ldist_reindex}
\ufsr \times (\ufss \sqcup \ufsu) \fto[\iso]{\upka}
(\ufsr \times \ufss) \sqcup (\ufsr \times \ufsu)
\end{equation}
induced by the universal properties of coproducts, sending $(i,k) \in \ufsr \times \ufss$ and $(i,p) \in \ufsr \times \ufsu$ in the domain to their counterparts in the codomain.  Using \cref{BA_mbrd_component} and the symmetry axiom \cref{smc_axioms} for the multiplicative braiding $\mbrd$, the component morphism of $\sldist$ is the identity morphism of the object
\begin{itemize}
\item $(\bigbtimes_{j \in \ufsm_i} a^i_j) \btimes (\bigbtimes_{l \in \ufsn_k} b^k_l)$ for each pair of indices $(i,k) \in \ufsr \times \ufss$ and
\item $(\bigbtimes_{j \in \ufsm_i} a^i_j) \btimes (\bigbtimes_{q \in \ufst_p} c^p_q)$ for each pair of indices $(i,p) \in \ufsr \times \ufsu$.
\end{itemize}
In each of the lower left diagram and the right diagram in \cref{biperm-axioms} for $(\sabrd,\srdist,\sldist)$, the two composite morphisms have
\begin{itemize}
\item the same reindexing functions by \cref{BA_mor_sum_re,BA_mor_prod_re,BA_abrd_reindex,BA_ldist_reindex}; and
\item identity component morphisms because each of $\sabrd$, $\srdist$, and $\sldist$ has identity component morphisms.
\end{itemize}
Thus, each diagram in \cref{biperm-axioms} commutes for the tuple
\[\big(\BsA, (\splus, \szero, \sabrd), (\stimes, \sunit, \smbrd), (\sladot, \srhodot), (\sldist, \srdist)\big),\]
proving that $\BsA$ is a bipermutative category \pcref{def:bipermcat}.
\end{definition}

\subsection*{Morphism Assignment of $\Bs$}

Having constructed the object assignment $\A \mapsto \BsA$ of $\Bs$, next we construct its morphism assignment.  

\begin{definition}\label{def:Bsf}
Suppose we are given a multiplicatively strong symmetric bimonoidal functor \cref{mult_strong} between bipermutative categories \pcref{def:bipermcat,def:sbfunctor}
\[\big(\sff,\ftwoplus,\fzeroplus,\ftwotimes,\fzerotimes\big) \cn \A \to \B,\]
with additive structure $\fplus = (\sff,\ftwoplus,\fzeroplus)$ and multiplicative structure $\ftimes = (\sff,\ftwotimes,\fzerotimes)$.  We define a \emph{strict} symmetric bimonoidal functor
\begin{equation}\label{Bsf}
\BsA \fto{\Bsf} \BsB
\end{equation}
as follows.
\begin{description}
\item[Objects] $\Bsf$ sends an object $\bba = \ang{\ang{a^i_j}_{j \in \ufsm_i}}_{i \in \ufsr}$ in $\BsA$, as defined in \cref{BA_object}, to the object
\begin{equation}\label{Bsf_object}
(\Bsf)\bba = \ang{\ang{\sff a^i_j}_{j \in \ufsm_i}}_{i \in \ufsr} \in \BsB
\end{equation}
by applying the functor $\sff$ to each alphabet $a^i_j$.  Thus, the object $(\Bsf)\bba$ has the same additive length as $\bba$, and its $i$-th monomial $\sff a^i = \ang{\sff a^i_j}_{j \in \ufsm_i}$ has the same multiplicative length as the $i$-th monomial $a^i = \ang{a^i_j}_{j \in \ufsm_i}$ of $\bba$.
\item[Morphisms] $\Bsf$ sends a morphism $(\uphi; g) \cn \bba \to \bbb$ in $\BsA$, as defined in \cref{BA_morphism}, to the morphism
\begin{equation}\label{Bsf_morphism}
\ang{\ang{\sff a^i_j}_{j \in \ufsm_i}}_{i \in \ufsr} = (\Bsf)\bba 
\fto[= (\uphi; (\Bsf)g)]{(\Bsf)(\uphi; g)} 
(\Bsf)\bbb = \ang{\ang{\sff b^k_l}_{l \in \ufsn_k}}_{k \in \ufss}
\end{equation}
in $\BsB$, with the same reindexing function $\uphi \cn \ufsr \to \ufss$.  The component morphism of $(\Bsf)g$ at an index $k \in \ufss$ is the following composite morphism in $\B$.
\begin{equation}\label{Bsf_mor_ii}
\begin{tikzpicture}[vcenter]
\def\v{-1.8} \def\h{3.5} \def\u{.8ex}
\draw[0cell=.85]
(0,0) node (a11) {\bigbplus_{i \in \uphiinv(k)} \bigbtimes_{j \in \ufsm_i} \sff a^i_j}
(a11)++(2*\h,0) node (a13) {\bigbtimes_{l \in \ufsn_k} \sff b^k_l}
(a11)++(0,\v) node (a21) {\bigbplus_{i \in \uphiinv(k)} \sff \Big(\bigbtimes_{j \in \ufsm_i} a^i_j\Big)}
(a21)++(\h+.3,0) node (a22) {\sff \Big(\bigbplus_{i \in \uphiinv(k)} \bigbtimes_{j \in \ufsm_i} a^i_j\Big)}
(a13)++(0,\v) node (a23) {\sff \Big(\bigbtimes_{l \in \ufsn_k} b^k_l\Big)}
;
\draw[1cell=.8]
(a11) edge[transform canvas={yshift=\u}] node {((\Bsf)g)^k} (a13)
(a11) edge node[swap] {\bigbplus_{i \in \uphiinv(k)} \ftwotimes} (a21)
(a21) edge[transform canvas={yshift=\u}] node {\ftwoplus} (a22)
(a22) edge[transform canvas={yshift=\u}] node {\sff g^k} (a23)
(a23) edge node[swap] {\ftwotimesinv} (a13)
;
\end{tikzpicture}
\end{equation}
\begin{itemize}
\item In the first arrow, $\ftwotimes$ denotes an iterate of the multiplicative monoidal constraint of $\sff$ \cref{mon_constraint}.  If $m_i = 0$, then it is defined as the multiplicative unit constraint $\fzerotimes \cn \bunit \to \sff \bunit$.
\item In the second arrow, $\ftwoplus$ denotes an iterate of the additive monoidal constraint of $\sff$.  If $\uphiinv(k) = \emptyset$, then it is defined as the additive unit constraint $\fzeroplus \cn \bzero \to \sff \bzero$.  If $|\uphiinv(k)| = 1$, then $\ftwoplus$ denotes the identity morphism.
\item In the third arrow, $g^k$ is the indicated component morphism of $(\uphi;g)$, as defined in \cref{BA_mor_ii}.
\item In the fourth arrow, $\ftwotimesinv$ denotes the inverse of an iterate of the multiplicative monoidal constraint of $\sff$ or, if $n_k = 0$, the inverse of the multiplicative unit constraint $\fzerotimes \cn \bunit \to \sff \bunit$.  
\end{itemize}
The arrow $\ftwotimesinv$ in \cref{Bsf_mor_ii} is the reason why we need to assume that the symmetric bimonoidal functor $\sff$ is multiplicatively strong.
\item[Functoriality] 
$\Bsf$ preserves identity morphisms because $\sff$ does.  Composition of morphisms \pcref{def:BAcomp} is preserved by $\Bsf$ by the uniqueness of coherence isomorphisms \cite[Theorem XI.1]{maclane} in the permutative category $\Aplus = (\A, \bplus, \bzero, \abrd)$ and the coherence axioms \cref{mfaxioms,ftwobraiding} for the additive structure $\fplus = (\sff,\ftwoplus,\fzeroplus)$.
\item[Additive structure]
By \cref{Bsf_object}, $\Bsf$ preserves the additive monoidal unit $\szero$ \cref{szero} and sums of objects \cref{splus}.  By \cref{Bsf_morphism}, $\Bsf$ preserves sums of morphisms \cref{BA_mor_sum} and the additive braiding $\sabrd$ \cref{BA_sabrd}.  Thus, defining the additive monoidal and unit constraints of $\Bsf$ to be identities, we have constructed a strict symmetric monoidal functor 
\[\BsAplus = (\BsA, \splus, \szero, \sabrd) \fto{\Bsf} \BsBplus = (\BsB, \splus, \szero, \sabrd)\]
between permutative categories.
\item[Multiplicative structure]
By \cref{Bsf_object}, $\Bsf$ preserves the multiplicative monoidal unit $\sunit$ \cref{sunit} and products of objects \cref{stimes}.  The morphism assignment $(\uphi; g) \mapsto (\Bsf)(\uphi; g)$ in \cref{Bsf_morphism} preserves
\begin{itemize}
\item products of morphisms \cref{BA_mor_prod} by \cref{Bsf_preserves_prod} below; and
\item the multiplicative braiding $\smbrd$ \cref{BA_smbrd} by the axiom \cref{ftwobraiding} for $\ftimes = (\sff,\ftwotimes)$, since $\ftwoplus$ in the diagram \cref{Bsf_mor_ii} is the identity when $\uphi$ is the swapping permutation $\tau$ in \cref{BA_mbrd_reindex}.
\end{itemize}
Thus, defining the multiplicative monoidal and unit constraints of $\Bsf$ to be identities, we have constructed a strict symmetric monoidal functor 
\[\BsAtimes = (\BsA, \stimes, \sunit, \smbrd) \fto{\Bsf} \BsBtimes = (\BsB, \stimes, \sunit, \smbrd)\]
between permutative categories.
\item[Coherence axioms]
Each of the two diagrams in \cref{sbf-axioms} commutes for $\Bsf$ because each arrow is the identity morphism.
\end{description}
This finishes the construction of the strict symmetric bimonoidal functor $\Bsf \cn \BsA \to \BsB$.
\end{definition}

\begin{lemma}\label{Bsf_preserves_prod}
In \cref{def:Bsf}, the morphism assignment $(\uphi; g) \mapsto (\Bsf)(\uphi; g)$ in \cref{Bsf_morphism} preserves products of morphisms \cref{BA_mor_prod}.
\end{lemma}

\begin{proof}
Using the notation in \cref{def:BA_mor_sumprod}, we need to show that the two morphisms
\begin{equation}\label{Bsf_pres_prod_ii}
\begin{tikzpicture}[vcenter]
\def\v{-.8}
\draw[0cell]
(0,0) node (a1) {(\Bsf)\bba \stimes (\Bsf)\bbc}
(a1)++(6,0) node (a2) {(\Bsf)\bbb \stimes (\Bsf)\bbd}
(a1)++(0,\v) node (b1) {(\Bsf)(\bba \stimes \bbc)}
(a2)++(0,\v) node (b2) {(\Bsf)(\bbb \stimes \bbd)}
;
\draw[1cell=.9]
(a1) edge[equal] (b1)
(a2) edge[equal] (b2)
(a1) edge node {(\Bsf)(\uphi; g) \stimes (\Bsf)(\uppsi; h)} (a2)
(b1) edge node {(\Bsf)((\uphi; g) \stimes (\uppsi; h))} (b2)
;
\end{tikzpicture}
\end{equation}
in $\BsA$ are equal.  By \cref{BA_mor_prod_re,Bsf_morphism}, each of these two morphisms has reindexing function given by the product $\uphi \times \uppsi$.  To verify that their component morphisms at a pair of indices $(k,y) \in \ufss \times \ufsw$ are equal, we assume that each index in the list below runs through the indicated ordered finite set.
\[\begin{aligned}
i & \in \uphiinv(k) \phantom{M} & j & \in \ufsm_i \phantom{M} & l & \in \ufsn_k \\
p & \in \uppsiinv(y) & q & \in \ufst_p & z & \in \ufsv_y
\end{aligned}\]
For example, $\bplus_i$ means $\bplus_{i \in \uphiinv(k)}$, and $\btimes_j$ means $\btimes_{j \in \ufsm_i}$.  By \cref{BA_mor_prod_ii,Bsf_mor_ii}, the $(k,y)$-component morphisms of the top and bottom morphisms in \cref{Bsf_pres_prod_ii} are the top and bottom boundary composites in the following diagram in $\B$.
\begin{equation}\label{Bsf_pres_prod_diag}
\begin{tikzpicture}[vcenter]
\def\e{.5} \def\g{1} \def\h{9} \def\u{1.2} \def\v{1.3}
\draw[0cell=.85]
(0,0) node (a11) {\bplus_i \bplus_p (\btimes_j \sff a^i_j) \btimes (\btimes_q \sff c^p_q)}
(a1)++(\h,0) node (a12) {(\btimes_l \sff b^k_l) \btimes (\btimes_z \sff d^y_z)}
(a11)++(\e,\v) node (a21) {(\bplus_i \btimes_j \sff a^i_j) \btimes (\bplus_p \btimes_q \sff c^p_q)}
(a12)++(-\e,\v) node (a22) {\sff(\btimes_l b^k_l) \btimes \sff(\btimes_z d^y_z)}
(a21)++(\g,\u) node (a31) {[\bplus_i \sff(\btimes_j a^i_j)] \btimes [\bplus_p \sff(\btimes_q c^p_q)]}
(a22)++(-\g,\u) node (a32) {\sff(\bplus_i \btimes_j a^i_j) \btimes \sff(\bplus_p \btimes_q c^p_q)} 
(a11)++(\e,-\v) node (b11) {\bplus_i \bplus_p \sff[(\btimes_j a^i_j) \btimes (\btimes_q c^p_q)]}
(a12)++(-\e,-\v) node (b12) {\sff[(\btimes_l b^k_l) \btimes (\btimes_z d^y_z)]}
(b11)++(\g,-\u) node (b21) {\sff[\bplus_i \bplus_p (\btimes_j a^i_j) \btimes (\btimes_q c^p_q)]}
(b12)++(-\g,-\u) node (b22) {\sff[(\bplus_i \btimes_j a^i_j) \btimes (\bplus_p \btimes_q c^p_q)]}
(a11)++(4,0) node (c) {\bplus_i \bplus_p \sff(\btimes_j a^i_j) \btimes \sff(\btimes_q c^p_q)}
;
\draw[1cell=.8]
(a11) edge node[pos=.1] {\deinv} (a21)
(a21) edge node[pos=.1] {(\bplus_i \ftwotimes) \btimes (\bplus_p \ftwotimes)} (a31)
(a31) edge node {\ftwoplus \btimes \ftwoplus} (a32)
(a32) edge node[pos=.9] {\sff g^k \btimes \sff h^y} (a22)
(a22) edge node[pos=.9] {\ftwotimesinv \btimes \ftwotimesinv} (a12)
(a11) edge node[swap,pos=.1] {\bplus_i \bplus_p \ftwotimes} (b11)
(b11) edge node[swap,pos=.1] {\ftwoplus} (b21)
(b21) edge node {\sff \deinv} (b22)
(b22) edge node[swap,pos=.9] {\sff(g^k \btimes h^y)} (b12)
(b12) edge node[swap,pos=.9] {\ftwotimesinv} (a12)
(a11) edge[bend left=10, shorten >=-1ex] node {\bplus_i \bplus_p \ftwotimes \btimes \ftwotimes} (c)
(a31) edge node {\de} (c)
(c) edge node[inner sep=0pt, pos=.4] {\bplus_i \bplus_p \ftwotimes} (b11)
(a32) edge[transform canvas={xshift=-4em}] node[swap,pos=.25] {\ftwotimes} (b22)
(a22) edge[transform canvas={xshift=-3em}] node[swap,pos=.5] {\ftwotimes} (b12)
;
\end{tikzpicture}
\end{equation}
Starting from the upper left quadrilateral and going counterclockwise, the five regions in the diagram \cref{Bsf_pres_prod_diag} commute for the following reasons:
\begin{itemize}
\item the naturality of Laplaza coherence isomorphisms \pcref{thm:laplaza};
\item the coherence axioms \cref{mfaxioms} for $\ftimes = (\sff,\ftwotimes,\fzerotimes)$ and the functoriality of $\bplus$; 
\item the coherence axioms \cref{sbf-axioms,sbf-axioms-left} for $\sff$;
\item the naturality of $\ftwotimes$; and
\item the coherence axioms \cref{mfaxioms} for $\ftimes$.
\end{itemize}
This proves that the two morphisms in \cref{Bsf_pres_prod_ii} are equal.
\end{proof}

\subsection*{Functoriality of $\Bs$}
The assignment $\sff \mapsto \Bsf$ in \cref{Bsf}  preserves identity morphisms by \cref{Bsf_object,Bsf_mor_ii}.  Recall from \cref{def:sbfunctor} that composition of symmetric bimonoidal functors are defined by separately composing the additive structures and the multiplicative structures.  Given another multiplicatively strong symmetric bimonoidal functor between bipermutative categories
\[(\sfh,\htwoplus,\hzeroplus,\htwotimes,\hzerotimes) \cn \B \to \C,\]
the diagram of strict symmetric bimonoidal functors
\[\begin{tikzpicture}
\def\h{1.8} \def\u{.6}
\draw[0cell]
(0,0) node (a1) {\BsA}
(a1)++(\h,0) node (a2) {\BsB}
(a2)++(\h,0) node (a3) {\BsC}
;
\draw[1cell=.9]
(a1) edge node {\Bsf} (a2)
(a2) edge node {\Bsh} (a3)
(a1) [rounded corners=2pt] |- ($(a2)+(-1,\u)$) -- node {\Bs(\sfh\sff)} ($(a2)+(1,\u)$) -| (a3)
;
\end{tikzpicture}\]
commutes on objects by \cref{Bsf_object}.  This diagram commutes on morphisms by \cref{Bsf_mor_ii} and the naturality of the additive monoidal and unit constraints, $\htwoplus$ and $\hzeroplus$.  In summary, we have constructed a functor
\begin{equation}\label{Bs_functor}
\bipermms \fto{\Bs} \bipermst.
\end{equation}
Its object assignment, $\A \mapsto \BsA$, is given in \cref{def:BAobject,def:BAmorphisms,def:BAcomp,def:BA_mor_sumprod,def:BA_braidings,def:BA_mzero_dist}.  Its morphism assignment, $\sff \mapsto \Bsf$, is given in \cref{def:Bsf}.  This proves the first assertion of \cref{thm:sbfstrict}.

\subsection*{The Adjunction $(\lefta,\righta)$}

To finish the proof of \cref{thm:sbfstrict}, next we construct the adjunction 
\[\begin{tikzpicture}
\def\b{22} \def\h{1.8}
\draw[0cell]
(0,0) node (a1) {\phantom{\A}}
(a1)++(-.15,0) node (a1') {\Bs\A}
(a1)++(\h,0) node (a2) {\A}
(a1)++(\h/2,0) node (a0) {\bot}
;
\draw[1cell=.9]
(a1) edge[bend left=\b] node {\lefta} (a2)
(a2) edge[bend left=\b] node {\righta} (a1)
;
\end{tikzpicture}\]
for a bipermutative category $\A$ \pcref{def:bipermcat}, starting with the left adjoint $\lefta$.  

\begin{definition}[Left Adjoint]\label{def:lefta}
We define an additively strict and multiplicatively strictly unital strong symmetric bimonoidal functor 
\[\BsA \fto{\lefta} \A\]
as follows.  Moreover, $\lefta$ is natural in strict symmetric bimonoidal functors; see \cref{lefta_nat_diagram}.
\begin{description}
\item[Objects] 
$\lefta$ sends an object $\bba = \ang{\ang{a^i_j}_{j \in \ufsm_i}}_{i \in \ufsr}$ in $\BsA$, as defined in \cref{BA_object}, to the object
\begin{equation}\label{lefta_object}
\lefta\bba = \bigbplus_{i \in \ufsr} \bigbtimes_{j \in \ufsm_i} a^i_j \in \A.
\end{equation}
We recall our convention that an empty product is the multiplicative monoidal unit $\bunit \in \A$, and an empty sum is the additive monoidal unit $\bzero \in \A$.  Thus, $\lefta$ sends the additive monoidal unit $\szero \in \BsA$ \cref{szero} to $\bzero \in \A$ and the multiplicative monoidal unit $\sunit \in \BsA$ to $\bunit \in \A$.  Moreover, $\lefta$ preserves sums of objects \cref{splus}.
\item[Morphisms]
Given a morphism $(\uphi; g) \cn \bba \to \bbb$ in $\BsA$, as defined in \cref{BA_morphism}, the morphism $\lefta(\uphi; g)$ is defined as the following composite in $\A$.
\begin{equation}\label{lefta_morphism}
\begin{tikzpicture}[vcenter]
\def\u{1} \def\t{.7ex}
\draw[0cell=.95]
(0,0) node (a1) {\lefta\bba = \bigbplus_{i \in \ufsr} \bigbtimes_{j \in \ufsm_i} a^i_j}
(a1)++(3.7,0) node (a2) {\bigbplus_{k \in \ufss} \bigbplus_{i \in \uphiinv(k)} \bigbtimes_{j \in \ufsm_i} a^i_j}
(a2)++(4,0) node (a3) {\bigbplus_{k \in \ufss} \bigbtimes_{l \in \ufsn_k} b^k_l = \lefta\bbb}
;
\draw[1cell=.85]
(a1) edge[transform canvas={yshift=\t}] node {\upsi} node[swap] {\iso} (a2)
(a2) edge[transform canvas={yshift=\t}] node {\bigbplus_{k \in \ufss} g^k} (a3)
(a1) [rounded corners=2pt] |- ($(a2)+(-1,\u)$) -- node {\lefta(\uphi; g)} ($(a2)+(1,\u)$) -| (a3)
;
\end{tikzpicture}
\end{equation}
\begin{itemize}
\item With respect to the reindexing function $\uphi \cn \ufsr \to \ufss$, the arrow $\upsi$ is the unique coherence isomorphism \cite[Theorem XI.1]{maclane} in the additive structure $\Aplus = (\A,\bplus,\bzero,\abrd)$.  It permutes the objects $\bigbtimes_{j \in \ufsm_i} a^i_j$ and inserts a copy of the additive monoidal unit $\bzero$ for each index $k \in \ufss$ with $\uphiinv(k) = \emptyset$.
\item In the second arrow, $g^k$ is the indicated component morphism of $g$, as defined in \cref{BA_mor_ii}.
\end{itemize}
\item[Functoriality]
The morphism assignment $(\uphi; g) \mapsto \lefta(\uphi;g)$ preserves identity morphisms and composition \cref{BA_mor_comp} by the naturality and uniqueness of coherence isomorphisms in the additive structure $\Aplus$.
\item[Additive structure]
The morphism assignment $(\uphi; g) \mapsto \lefta(\uphi;g)$ preserves sums of morphisms \cref{BA_mor_sum} and the additive braiding $\sabrd$ \cref{BA_sabrd} by the functoriality of the sum $\bplus$ of $\A$.  Thus, defining the additive monoidal and unit constraints of $\lefta$ to be identities, we have constructed a strict symmetric monoidal functor between permutative categories
\[\BsAplus = (\BsA, \splus, \szero, \sabrd) 
\fto{\lefta} \Aplus = (\A,\bplus,\bzero,\abrd).\]
\item[Multiplicative structure]
We define the multiplicative unit constraint of $\lefta$ to be the identity morphism.  It is well defined because $\lefta\sunit = \bunit$.  

The multiplicative monoidal constraint $\leftatwotimes$ for two objects $\bba = \ang{\ang{a^i_j}_{j \in \ufsm_i}}_{i \in \ufsr}$ and $\bbb = \ang{\ang{b^k_l}_{l \in \ufsn_k}}_{k \in \ufss}$ in $\BsA$ is defined as the following unique Laplaza coherence isomorphism in $\A$ \pcref{thm:laplaza}, where the product $\bba \stimes \bbb = \ang{\ang{a^i b^k}_{k \in \ufss}}_{i \in \ufsr}$ is defined in \cref{stimes}.
\begin{equation}\label{leftatwotimes}
\begin{tikzpicture}[vcenter]
\def\v{-1}
\draw[0cell=.95]
(0,0) node (a11) {\lefta\bba \btimes \lefta\bbb}
(a11)++(5.2,0) node (a12) {\lefta(\bba \stimes \bbb)}
(a11)++(0,\v) node (a21) {\Big(\bigbplus_{i \in \ufsr} \bigbtimes_{j \in \ufsm_i} a^i_j\Big) \btimes \Big(\bigbplus_{k \in \ufss} \bigbtimes_{l \in \ufsn_k} b^k_l\Big)}
(a12)++(0,\v) node (a22) {\bigbplus_{i \in \ufsr} \bigbplus_{k \in \ufss} \Big(\bigbtimes_{j \in \ufsm_i} a^i_j\Big) \btimes \Big(\bigbtimes_{l \in \ufsn_k} b^k_l\Big)}
;
\draw[1cell=.85]
(a11) edge[equal] (a21)
(a12) edge[equal] (a22)
(a11) edge node {\leftatwotimes} (a12)
(a21) edge[transform canvas={yshift=.5ex}] node {\de} node[swap] {\iso} (a22)
;
\end{tikzpicture}
\end{equation}
\begin{itemize}
\item If $\bba, \bbb \neq \szero$, then the isomorphism $\de$ consists of an iterate of the right distributivity $\rdist$ of $\A$, which is the identity, and, for each $i \in \ufsr$, an iterate of the left distributivity $\ldist$ \cref{larhordist}.
\item If $\bba = \szero$, then 
\[\bba \stimes \bbb = \szero \andspace \lefta\bba = \lefta(\bba \stimes \bbb) = \bzero.\]  
In this case, $\de$ is the left multiplicative zero $\ladot$ of $\A$, which is the identity.
\item If $\bbb = \szero$, then 
\[\bba \stimes \bbb = \szero \andspace \lefta\bbb = \lefta(\bba \stimes \bbb) = \bzero.\]  
In this case, $\de$ is the right multiplicative zero $\rhodot$ of $\A$, which is the identity.
\end{itemize}

The naturality of $\leftatwotimes$ is verified in \cref{ltwotimes_nat} below.  Each of the four axioms of a symmetric monoidal functor in \cref{mfaxioms,ftwobraiding} holds for $(\lefta,\leftatwotimes,1_{\bunit})$ by the uniqueness of Laplaza coherence isomorphisms.  Thus, we have constructed a strictly unital strong symmetric monoidal functor between permutative categories
\[\BsAtimes = (\BsA, \stimes, \sunit, \smbrd) 
\fto{(\lefta,\leftatwotimes,1_{\bunit})} \Atimes = (\A,\btimes,\bunit,\mbrd).\]
\item[Coherence axioms]
The left diagram in \cref{sbf-axioms} commutes for $\lefta$ because each composite is equal to $\rhodot$.  The right diagram in \cref{sbf-axioms} commutes for $\lefta$ by the uniqueness of Laplaza coherence isomorphisms.
\item[Naturality]
The naturality of $\lefta$ in strict symmetric bimonoidal functors $\sff \cn \A \to \B$ means that the diagram of symmetric bimonoidal functors
\begin{equation}\label{lefta_nat_diagram}
\begin{tikzpicture}[vcenter]
\def\v{-1.3}
\draw[0cell]
(0,0) node (a11) {\BsA}
(a11)++(1.8,0) node (a12) {\A}
(a11)++(0,\v) node (a21) {\BsB}
(a12)++(0,\v) node (a22) {\B}
;
\draw[1cell=.9]
(a11) edge node {\lefta} (a12)
(a12) edge node {\sff} (a22)
(a11) edge node[swap] {\Bsf} (a21)
(a21) edge node {\lefta} (a22)
;
\end{tikzpicture}
\end{equation}
commutes.
\begin{itemize}
\item The diagram \cref{lefta_nat_diagram} commutes on objects by \cref{Bsf_object}, \cref{lefta_object}, and the strictness of $\sff$.
\item It commutes on morphisms by \cref{Bsf_morphism}, \cref{lefta_morphism}, and the strictness of $\sff$.
\item To see that the two composites in the diagram \cref{lefta_nat_diagram} have the same additive and multiplicative structures, we use the strictness of the symmetric bimonoidal functor $\Bsf$ \pcref{def:Bsf}, the definition \cref{leftatwotimes} of $\leftatwotimes$, and the coherence axioms for $\sff$ in \cref{sbf-axioms,sbf-axioms-left}.
\end{itemize}
\end{description}
This completes the construction of the additively strict and multiplicatively strictly unital strong symmetric bimonoidal functor $\lefta \cn \BsA \to \A$ that is natural in strict symmetric bimonoidal functors.
\end{definition}

\begin{lemma}\label{ltwotimes_nat}
In \cref{leftatwotimes}, $\leftatwotimes$ is a natural transformation.
\end{lemma}

\begin{proof}
Using the notation in \cref{BA_mor_prod}, we need to show that the following diagram in $\A$ commutes.
\begin{equation}\label{twotimes_nat_diag}
\begin{tikzpicture}[vcenter]
\def\v{-1.3}
\draw[0cell]
(0,0) node (a11) {\lefta\bba \btimes \lefta\bbc}
(a11)++(2.5,0) node (a12) {\lefta(\bba \stimes \bbc)}
(a11)++(0,\v) node (a21) {\lefta\bbb \btimes \lefta\bbd}
(a12)++(0,\v) node (a22) {\lefta(\bbb \stimes \bbd)}
;
\draw[1cell=.9]
(a11) edge node {\leftatwotimes} (a12)
(a12) edge node {\lefta((\uphi;g) \stimes (\uppsi;h))} (a22)
(a11) edge node[swap] {\lefta(\uphi;g) \btimes \lefta(\uppsi;h)} (a21)
(a21) edge node {\leftatwotimes} (a22)
;
\end{tikzpicture}
\end{equation}
We assume that each index in the list below runs through the indicated ordered finite set.
\[\begin{aligned}
i & \in \ufsr & j & \in \ufsm_i & p & \in \ufsu & q & \in \ufst_p\\
k & \in \ufss & l & \in \ufsn_k & y & \in \ufsw & z & \in \ufsv_y\\
i' & \in \uphiinv(k) &&& p' & \in \uppsiinv(y) && 
\end{aligned}\]
By \cref{BA_mor_prod,lefta_morphism,leftatwotimes}, the diagram \cref{twotimes_nat_diag} is given by the boundary diagram below.
\[\begin{tikzpicture}
\def\v{-1.3}
\draw[0cell=.9]
(0,0) node (a11) {(\bplus_i \btimes_j a^i_j) \btimes (\bplus_p \btimes_q c^p_q)}
(a11)++(5.6,0) node (a12) {\bplus_i \bplus_p (\btimes_j a^i_j) \btimes (\btimes_q c^p_q)}
(a12)++(0,\v) node (a22) {\bplus_k \bplus_y \bplus_{i'} \bplus_{p'} (\btimes_j a^{i'}_j) \btimes (\btimes_q c^{p'}_q)}
(a11)++(0,2*\v) node (a31) {(\bplus_k \bplus_{i'} \btimes_j a^{i'}_j) \btimes (\bplus_y \bplus_{p'} \btimes_q c^{p'}_q)}
(a22)++(0,\v) node (a32) {\bplus_k \bplus_y (\bplus_{i'} \btimes_j a^{i'}_j) \btimes (\bplus_{p'} \btimes_q c^{p'}_q)}
(a31)++(0,\v) node (a41) {(\bplus_k \btimes_l b^k_l) \btimes (\bplus_y \btimes_z d^y_z)}
(a32)++(0,\v) node (a42) {\bplus_k \bplus_y (\btimes_l b^k_l) \btimes (\btimes_z d^y_z)}
;
\draw[1cell=.85]
(a11) edge node[swap] {\si \btimes \si} (a31)
(a31) edge node[swap] {(\bplus_k g^k) \btimes (\bplus_y h^y)} (a41)
(a41) edge node {\de} (a42)
(a11) edge node {\de} (a12)
(a12) edge node {\si} (a22)
(a22) edge node {\bplus_k \bplus_y \deinv} (a32)
(a32) edge node {\bplus_k \bplus_y  (g^k \btimes h^y)} (a42)
(a31) edge node {\de} (a32)
;
\end{tikzpicture}\]
In the previous diagram, the top and bottom rectangles commute by, respectively, the uniqueness and naturality of Laplaza coherence isomorphisms \pcref{thm:laplaza}.
\end{proof}

Next, we define the right adjoint of $\lefta$.

\begin{definition}[Right Adjoint]\label{def:righta}
Given a bipermutative category $\A$, we define a multiplicatively strong symmetric bimonoidal functor 
\[\A \fto{\righta} \BsA\]
as follows.  Moreover, $\righta$ is natural in multiplicatively strong symmetric bimonoidal functors; see \cref{righta_nat_diagram}.
\begin{description}
\item[Objects]
Using the terminology in \cref{BA_object}, $\righta$ sends an object $a \in \A$ to the object
\begin{equation}\label{righta_object}
\righta a = \ang{\ang{a}} \in \BsA
\end{equation}
with additive length 1.  Its only monomial has multiplicative length 1 and alphabet $a$.  To simplify the notation below, we abbreviate $\ang{\ang{a}}$ to $\ang{a}$.
\item[Morphisms]
Using the terminology in \cref{def:BAmorphisms}, $\righta$ sends a morphism $g \cn a \to b$ in $\A$ to the morphism
\begin{equation}\label{righta_morphism}
\righta a = \ang{a} \fto{(1_{\ufsone}; g)} \righta b = \ang{b} \inspace \BsA
\end{equation}
with identity reindexing function $1 \cn \ufsone \to \ufsone$.  Its only component morphism is $g$.  This defines a functor $\righta \cn \A \to \BsA$ because, in the diagram \cref{BA_comp_hg}, the coherence isomorphism $\si_{\uphi, \uppsi;p}$ is the identity when both $\uphi$ and $\uppsi$ are the identity morphism on $\ufsone$.
\item[Additive structure]
With the additive monoidal unit $\szero$ in \cref{szero}, the additive unit constraint of $\righta$ is the morphism
\begin{equation}\label{rightazeroplus}
\szero \fto{\rightazeroplus = (\uphi_0; 1_\bzero)} 
\righta\bzero = \ang{\bzero} \inspace \BsA
\end{equation}
with reindexing function $\uphi_0 \cn \ufs{0} = \emptyset \to \ufsone$ given by the unique function to $\ufsone$.  Recalling our convention that an empty sum is the additive monoidal unit $\bzero \in \A$, the only component morphism of $\rightazeroplus$ is the identity $1_\bzero \cn \bzero \to \bzero$.

With the sum of objects in \cref{splus}, the additive monoidal constraint of $\righta$ for two objects $a,b \in \A$ is the morphism
\begin{equation}\label{rightatwoplus}
\righta a \splus \righta b = \ang{a} \splus \ang{b} = \ang{\ang{a}, \ang{b}}
 \fto{\rightatwoplus = (\uphi_2; 1_{a \bplus b})} 
\righta(a \bplus b) = \ang{a \bplus b}
\end{equation}
in $\BsA$ with reindexing function $\uphi_2 \cn \ufs{2} \to \ufsone$ given by the unique function to $\ufsone$.  Its only component morphism is the identity morphism $1_{a \bplus b}$ of the sum $a \bplus b$.  The naturality of $\rightatwoplus$ follows from \cref{BA_comp_hg,BA_mor_sum,righta_morphism}.

This defines a symmetric monoidal functor between permutative categories
\[\Aplus = (\A,\bplus,\bzero,\abrd) \fto{(\righta,\rightatwoplus,\rightazeroplus)} 
\BsAplus = (\BsA,\splus,\szero,\sabrd).\]
In each of the four diagrams in \cref{mfaxioms,ftwobraiding} for $(\righta,\rightatwoplus,\rightazeroplus)$, each composite reindexing function is the unique function to $\ufsone$.  In \cref{mfaxioms}, each composite component morphism is the identity.  In \cref{ftwobraiding}, each composite component morphism is the additive braiding $\abrd \cn a \bplus b \to b \bplus a$.
\item[Multiplicative structure]
With the multiplicative monoidal unit $\sunit$ in \cref{sunit}, the multiplicative unit constraint of $\righta$ is the isomorphism
\begin{equation}\label{rightazerotimes}
\sunit = \ang{\emptyset}
\fto[\iso]{\rightazerotimes = (1_{\ufsone}; 1_{\bunit})} 
\righta\bunit = \ang{\bunit} \inspace \BsA
\end{equation}
with identity reindexing function $1 \cn \ufsone \to \ufsone$.  Recalling our convention that an empty product is the multiplicative monoidal unit $\bunit \in \A$, the only component morphism of $\rightazerotimes$ is the identity of $\bunit$.  Note that $\rightazerotimes$ is not the identity morphism because the only monomial in its domain, respectively, codomain, has multiplicative length 0, respectively, 1.

With the product of objects in \cref{stimes}, the multiplicative monoidal constraint of $\righta$ for two objects $a,b \in \A$ is the isomorphism
\begin{equation}\label{rightatwotimes}
\righta a \stimes \righta b = \ang{a} \stimes \ang{b} = \ang{\ang{a,b}} 
\fto[\iso]{\rightatwotimes = (1_{\ufsone}; 1_{a \btimes b})}
\righta(a \btimes b) = \ang{a \btimes b}
\end{equation}
with identity reindexing function $1 \cn \ufsone \to \ufsone$.  The domain of $\rightatwotimes$, $\ang{\ang{a,b}}$, has additive length 1, and its only monomial, $\ang{a,b}$, has multiplicative length 2.  The only component morphism of $\rightatwotimes$ is the identity morphism $1_{a \btimes b}$ of the product $a \btimes b$.  The naturality of $\rightatwotimes$ follows from \cref{BA_comp_hg}, \cref{righta_morphism}, and the fact that, in the diagram \cref{BA_mor_prod_ii}, the Laplaza coherence isomorphism $\deinv_{\uphi,\uppsi; k,y}$ is the identity when both $\uphi$ and $\uppsi$ are the identity morphism on $\ufsone$.  Note that $\rightatwotimes$ is not the identity morphism because the only monomial in its domain, respectively, codomain, has multiplicative length 2, respectively, 1.

This defines a strong symmetric monoidal functor between permutative categories
\[\Atimes = (\A,\btimes,\bunit,\mbrd) \fto{(\righta,\rightatwotimes,\rightazerotimes)} 
\BsAtimes = (\BsA,\stimes,\sunit,\smbrd).\]
In each of the four diagrams in \cref{mfaxioms,ftwobraiding} for $(\righta,\rightatwotimes,\rightazerotimes)$, each composite reindexing function is the unique function to $\ufsone$.  In \cref{mfaxioms}, each composite component morphism is the identity.  In \cref{ftwobraiding}, each composite component morphism is the multiplicative braiding $\mbrd \cn a \btimes b \to b \btimes a$.
\item[Coherence axioms]
The left diagram in \cref{sbf-axioms} commutes for $\righta$ because the multiplicative monoidal constraint
\[\righta a \stimes \righta \bzero = \ang{\ang{a,\bzero}} 
\fto[\iso]{\rightatwotimes = (1_{\ufsone}; 1_{a \btimes \bzero})}
\righta(a \btimes \bzero) = \ang{a \btimes \bzero} = \ang{\bzero}\]
has component morphism given by the identity of $a \btimes \bzero = \bzero$.  In the right diagram in \cref{sbf-axioms} for $\righta$, each composite reindexing function is the unique function to $\ufsone$, and each composite component morphism is the identity of the object 
\[(a \bplus b) \btimes c = (a \btimes c) \bplus (b \btimes c).\]
\item[Naturality]
The naturality of $\righta$ in multiplicatively strong symmetric bimonoidal functors $\sff \cn \A \to \B$ means that the diagram of symmetric bimonoidal functors
\begin{equation}\label{righta_nat_diagram}
\begin{tikzpicture}[vcenter]
\def\v{-1.3}
\draw[0cell]
(0,0) node (a11) {\A}
(a11)++(1.8,0) node (a12) {\BsA}
(a11)++(0,\v) node (a21) {\B}
(a12)++(0,\v) node (a22) {\BsB}
;
\draw[1cell=.9]
(a11) edge node {\righta} (a12)
(a12) edge node {\Bsf} (a22)
(a11) edge node[swap] {\sff} (a21)
(a21) edge node {\righta} (a22)
;
\end{tikzpicture}
\end{equation}
commutes.  Recall from \cref{def:Bsf} that $\Bsf$ is only defined when the symmetric bimonoidal functor $\sff$ is multiplicatively strong because the diagram \cref{Bsf_mor_ii} involves the arrow $\ftwotimesinv$.
\begin{itemize}
\item The diagram \cref{righta_nat_diagram} commutes on objects by \cref{Bsf_object,righta_object}.
\item It commutes on morphisms by \cref{Bsf_morphism,righta_morphism}.
\item To see that the two composites in the diagram \cref{righta_nat_diagram} have the same additive and multiplicative structures, the keys are the diagram \cref{Bsf_mor_ii} that defines the morphism assignment of $\Bsf$ and the fact that $\Bsf$ is strict.  For each of the two composites in \cref{righta_nat_diagram} and for the additive or multiplicative structures, the unit or monoidal constraints have component morphisms given by the corresponding constraints for $(\sff, \ftwoplus, \fzeroplus, \ftwotimes, \fzerotimes)$, while the reindexing function is the unique function to $\ufsone$.
\end{itemize}
\end{description}
We have constructed the multiplicatively strong symmetric bimonoidal functor $\righta \cn \A \to \BsA$ that is natural in multiplicatively strong symmetric bimonoidal functors.
\end{definition}

We can now finish the proof of our main theorem.

\begin{proof}[Proof of \cref{thm:sbfstrict}]
The functors $\Bs$ \cref{Bs_functor}, $\lefta \cn \BsA \to \A$ \pcref{def:lefta}, and $\righta \cn \A \to \BsA$ \pcref{def:righta} have already been constructed.  It remains to construct the unit and the identity counit for the bimonoidal adjunction $(\lefta,\righta)$, and verify the triangle identities in \cite[IV.1 Theorem 2(v)]{maclane}.

\parhead{Counit}.  For a bipermutative category $\A$, the composite
\begin{equation}\label{rl_id}
\A \fto{\righta} \BsA \fto{\lefta} \A
\end{equation}
is the identity symmetric bimonoidal functor \pcref{def:sbfunctor} for the following reasons. 
\begin{itemize}
\item $\lefta\righta$ is the identity on objects by \cref{lefta_object,righta_object}, and the identity on morphisms by \cref{lefta_morphism,righta_morphism}. 
\item The additive unit constraint, additive monoidal constraint, and multiplicative unit constraint of $\lefta\righta$ are identities because these constraints are identities for $\lefta$, along with the definitions \cref{monfunctor_comp,lefta_morphism,rightazeroplus,rightatwoplus,rightazerotimes}. 
\item The multiplicative monoidal constraint of $\lefta\righta$ for objects $a,b \in \A$ is the composite
\[a \btimes b = \lefta\righta a \btimes \lefta\righta b 
\fto{\leftatwotimes} \lefta(\righta a \stimes \righta b) = a \btimes b
\fto{\lefta(\rightatwotimes)} \lefta\righta(a \btimes b).\]
By \cref{leftatwotimes} and the uniqueness of Laplaza coherence isomorphisms \pcref{thm:laplaza}, the first arrow $\leftatwotimes$ is the identity.  The second arrow $\lefta(\rightatwotimes)$ is the identity by \cref{lefta_morphism,rightatwotimes}.
\end{itemize}
We define the counit 
\begin{equation}\label{counita}
\begin{tikzpicture}[vcenter]
\def\b{28} \def\h{1.8}
\draw[0cell]
(0,0) node (a1) {\A}
(a1)++(\h,0) node (a2) {\A}
;
\draw[1cell=.9]
(a1) edge[bend left=\b] node {\lefta\righta} (a2)
(a1) edge[bend right=\b] node[swap] {1_{\A}} (a2)
;
\draw[2cell=.9]
node[between=a1 and a2 at .45, rotate=-90, 2labelmed={above,\counita}] {\Rightarrow}
;
\end{tikzpicture}
\end{equation}
to be the identity bimonoidal natural transformation \pcref{def:bimonnat}.

\parhead{Unit}. We define the unit to be the bimonoidal natural transformation
\[\begin{tikzpicture}
\def\b{28} \def\h{1.8}
\draw[0cell]
(0,0) node (a1) {\phantom{\A}}
(a1)++(\h,0) node (a2) {\phantom{\A}}
(a1)++(-.15,0) node (a1') {\BsA}
(a2)++(.15,0) node (a2') {\BsA}
;
\draw[1cell=.9]
(a1) edge[bend left=\b] node {1_{\BsA}} (a2)
(a1) edge[bend right=\b] node[swap] {\righta\lefta} (a2)
;
\draw[2cell=.9]
node[between=a1 and a2 at .45, rotate=-90, 2labelmed={above,\unita}] {\Rightarrow}
;
\end{tikzpicture}\]
whose component at an object $\bba = \ang{\ang{a^i_j}_{j \in \ufsm_i}}_{i \in \ufsr} \in \BsA$ is the morphism
\begin{equation}\label{unita}
\bba \fto{\unita_\bba = (\uphi_r; 1)} 
\righta\lefta\bba = \ang{\lefta\bba} 
= \bang{\textstyle\bigbplus_{i \in \ufsr} \bigbtimes_{j \in \ufsm_i} a^i_j}
\end{equation}
in $\BsA$ \pcref{def:BAmorphisms}.  The reindexing function of $\unita_\bba$ is the unique function $\uphi_r \cn \ufsr \to \ufsone$.  The only component morphism of $\unita_\bba$ is the identity morphism of the object $\bigbplus_{i \in \ufsr} \bigbtimes_{j \in \ufsm_i} a^i_j$.

\parhead{Naturality}.  The naturality of $\unita \cn 1_{\BsA} \to \righta\lefta$ means that, for each morphism $(\uphi; g) \cn \bba \to \bbb$ in $\BsA$, as defined in \cref{BA_morphism}, the diagram
\begin{equation}\label{unita_natural}
\begin{tikzpicture}[vcenter]
\def\v{-1.3}
\draw[0cell]
(0,0) node (a11) {\bba}
(a11)++(4.6,0) node (a12) {\righta\lefta\bba = \bang{\textstyle\bigbplus_{i \in \ufsr} \bigbtimes_{j \in \ufsm_i} a^i_j}}
(a11)++(0,\v) node (a21) {\bbb}
(a12)++(0,\v) node (a22) {\phantom{\righta\lefta\bbb = \bang{\textstyle\bigbplus_{k \in \ufss} \bigbtimes_{l \in \ufsn_k} b^k_l}}}
(a22)++(0,-.04) node (a22') {\righta\lefta\bbb = \bang{\textstyle\bigbplus_{k \in \ufss} \bigbtimes_{l \in \ufsn_k} b^k_l}}
;
\draw[1cell=.9]
(a11) edge node {\unita_\bba = (\uphi_r; 1)} (a12)
(a12) edge[transform canvas={xshift=-4em}, shorten <=-1ex, shorten >=-.5ex] node {(1_{\ufsone}; \lefta(\uphi;g))} (a22)
(a11) edge node[swap] {(\uphi; g)} (a21)
(a21) edge node {\unita_\bbb = (\uphi_s; 1)} (a22)
;
\end{tikzpicture}
\end{equation}
in $\BsA$ commutes. 
\begin{itemize}
\item By \cref{BA_comp_re,righta_morphism}, each of the two composites in the diagram \cref{unita_natural} has reindexing function given by the unique function $\ufsr \to \ufsone$. 
\item By \cref{BA_comp_hg}, each of the two composites in \cref{unita_natural} has component morphism given by $\lefta(\uphi;g)$ in \cref{lefta_morphism}.
\end{itemize}

\parhead{Bimonoidality}.  The natural transformation $\unita \cn 1_{\BsA} \to \righta\lefta$ is bimonoidal \pcref{def:bimonnat} for the following reasons.
\begin{itemize}
\item In \cref{bimonnat_axioms} for $\unita$, the diagrams for the additive unit constraints, additive monoidal constraints, and multiplicative unit constraints commute because these constraints are identities for $\lefta$, along with the definitions \cref{monfunctor_comp,BA_mor_comp,BA_mor_sum,righta_morphism,rightazeroplus,rightatwoplus,rightazerotimes,unita}. 
\item The multiplicative monoidal constraint diagram in \cref{bimonnat_axioms} for $\unita \cn 1_{\BsA} \to \righta\lefta$ is the following diagram in $\BsA$ for objects $\bba = \ang{\ang{a^i_j}_{j \in \ufsm_i}}_{i \in \ufsr}$ and $\bbb = \ang{\ang{b^k_l}_{l \in \ufsn_k}}_{k \in \ufss}$.
\begin{equation}\label{unita_mmon}
\begin{tikzpicture}[vcenter]
\def\v{-1.4}
\draw[0cell]
(0,0) node (a11) {\bba \stimes \bbb}
(a11)++(3.5,0) node (a12) {\righta\lefta(\bba \stimes \bbb)}
(a11)++(0,\v) node (a21) {\righta\lefta\bba \stimes \righta\lefta\bbb}
(a12)++(0,\v) node (a22) {\righta(\lefta\bba \btimes \lefta\bbb)}
;
\draw[1cell=.9]
(a11) edge node {\unita_{\bba {\scalebox{.7}{$\stimes$}} \bbb}} (a12)
(a11) edge node[swap] {\unita_\bba \stimes \unita_\bbb} (a21)
(a21) edge node {\rightatwotimes} (a22)
(a22) edge node[swap] {\righta(\leftatwotimes)} (a12)
;
\end{tikzpicture}
\end{equation}
By the definitions \cref{stimes,BA_mor_prod,leftatwotimes,righta_morphism,rightatwotimes,unita}, the diagram \cref{unita_mmon} is given more explicitly by the diagram \cref{unita_mmon_ii} below, where we use the following abbreviations. 
\[a^i_\crdot = \dbigbtimes_{j \in \ufsm_i} a^i_j \qquad \dbigbplus_i = \dbigbplus_{i \in \ufsr} \qquad
b^k_\crdot = \dbigbtimes_{l \in \ufsn_k} b^k_l \qquad \dbigbplus_k = \dbigbplus_{k \in \ufss}\]
\begin{equation}\label{unita_mmon_ii}
\begin{tikzpicture}[vcenter]
\def\v{-1.5}
\draw[0cell=.9]
(0,0) node (a11) {\ang{\ang{a^i b^k}_{k \in \ufss}}_{i \in \ufsr}}
(a11)++(4.7,0) node (a12) {\bang{\dbigbplus_i \,\dbigbplus_k (a^i_\crdot \btimes b^k_\crdot)}}
(a11)++(0,\v) node (a21) {\bang{\bang{\dbigbplus_i a^i_\crdot \scs \dbigbplus_k b^k_\crdot}}}
(a12)++(0,\v) node (a22) {\bang{\big(\dbigbplus_i a^i_\crdot\big) \btimes \big(\dbigbplus_k b^k_\crdot\big)}}
;
\draw[1cell=.85]
(a11) edge node {(\uphi_{rs}; 1)} (a12)
(a11) edge node[swap] {(\uphi_{rs}; \deinv)} (a21)
(a21) edge node {(1_{\ufsone}; 1)} (a22)
(a22) edge node[swap] {(1_{\ufsone}; \de)} (a12)
;
\end{tikzpicture}
\end{equation}
By \cref{BA_mor_comp}, the diagram \cref{unita_mmon_ii} commutes.
\end{itemize}
This finishes the construction of the bimonoidal natural transformation $\unita \cn 1_{\BsA} \to \righta\lefta$.  Having defined the unit $\unita$ and the counit $\counita$ for $(\lefta,\righta)$, next we verify the triangle identities for an adjunction.

\parhead{Left triangle identity}.
For each object $\bba = \ang{\ang{a^i_j}_{j \in \ufsm_i}}_{i \in \ufsr} \in \BsA$, we need to show that the composite morphism
\begin{equation}\label{left_triangle}
\bigbplus_{i \in \ufsr} \bigbtimes_{j \in \ufsm_i} a^i_j 
= \lefta\bba \fto{\lefta\unita_\bba} \lefta\righta\lefta\bba 
\fto{\counita_{\lefta\bba}} \lefta\bba
\end{equation}
in $\A$ is the identity. 
\begin{itemize}
\item The second morphism $\counita_{\lefta\bba}$ is the identity by the definition \cref{counita} of the counit $\counita$. 
\item The first morphism 
\[\lefta\unita_\bba = \lefta(\uphi_r; 1)\] 
is the identity by \cref{unita}, since the coherence isomorphism $\upsi$ in \cref{lefta_morphism} is the identity for the reindexing function $\uphi_r \cn \ufsr \to \ufsone$.
\end{itemize}  
Thus, the composite in the diagram \cref{left_triangle} is the identity morphism.

\parhead{Right triangle identity}.
For each object $a \in \A$, we need to show that the composite morphism
\begin{equation}\label{right_triangle}
\ang{a} = \righta a \fto{\unita_{\righta a}} 
\righta\lefta\righta a \fto{\righta\counita_a} \righta a
\end{equation}
in $\BsA$ is the identity. 
\begin{itemize}
\item The second morphism $\righta\counita_a$ is the identity by the definition \cref{counita} of the counit $\counita$ and the functoriality of $\righta$. 
\item The first morphism 
\[\unita_{\righta a} = \unita_{\ang{a}} = (\uphi_1; 1) \cn \ang{a} \to \ang{a}\]
is the identity, since $\uphi_1 \cn \ufsone \to \ufsone$ is the identity function. 
\end{itemize} 
Thus, the composite in the diagram \cref{right_triangle} is the identity morphism.

In summary, for each bipermutative category $\A$, the quadruple $(\lefta,\righta,\unita,\counita)$ is a bimonoidal adjunction with left adjoint $\lefta$ \pcref{def:lefta}, right adjoint $\righta$ \pcref{def:righta}, unit $\unita$ \cref{unita}, and counit $\counita$ \cref{counita}.  This completes the proof of \cref{thm:sbfstrict}.
\end{proof}

\begin{remark}[Non-invertibility of the Unit]\label{rk:unita_notiso}
The unit $\unita \cn 1_{\BsA} \to \righta\lefta$ is not a natural isomorphism because a typical component $\unita_{\bba}$ \cref{unita} goes from an object $\bba$ with additive length $r \geq 0$ to an object $\righta\lefta\bba$ with additive length $1$.  Thus, the adjunction $(\lefta,\righta,\unita,\counita)$ in \cref{thm:sbfstrict} is generally not an adjoint equivalence.  This situation is similar to May's \cref{may4.3}, where the unit $1_{\Pst\A} \to \righta\lefta$ is not an isomorphism, as we mention in the paragraph of \cref{kr_weq}.
\end{remark}

We end this section with the following remark that illustrates some subtle differences between May's \cref{may4.3} in the permutative case and \cref{thm:sbfstrict} in the bipermutative case.  

\begin{remark}[Homotopy Theories]\label{rk:Bs_adjoint}
For the permutative case, the strictification functor $\Pst$ in May's \cref{may4.3} is the left adjoint of an adjoint equivalence of homotopy theories \cite[Def.\ 1.8]{gjo-extending}
\[\begin{tikzpicture}
\def\h{1.8} 
\draw[0cell]
(0,0) node (a1) {\phantom{X}}
(a1)++(\h,0) node (a2) {\phantom{X}}
(a1)++(-.25,0) node (a1') {\perm}
(a2)++(.37,0) node (a2') {\permst} 
(a1)++(\h/2,.02) node (a) {\bot}
;
\draw[1cell=.9]
(a1) edge[bend left=20] node {\Pst} (a2)
(a2) edge[bend left=20] node {\PIn} (a1)
;
\end{tikzpicture}\]
with $\PIn$ denoting the subcategory inclusion, unit $\righta \cn 1 \to \PIn\Pst$, and counit $\lefta \cn \Pst\PIn \to 1$.   On each side, the weak equivalences are defined as the morphisms that are sent to weak equivalences of spectra by May's $K$-theory machine \cite{may,may-groupcompletion,may-permutative}.

For the bipermutative case, no such equivalences of homotopy theories are known to exist.  The functor $\Bs$ in \cref{thm:sbfstrict} is \emph{almost}, but not quite, a part of an adjunction
\[\begin{tikzpicture}
\def\h{1.6} \def\m{.6}
\draw[0cell]
(0,0) node (a1) {\phantom{X}}
(a1)++(\h,0) node (a2) {\phantom{X}}
(a1)++(-\m,0) node (a1') {\bipermms}
(a2)++(.55,0) node (a2') {\bipermst} 
;
\draw[1cell=.9]
(a1) edge[transform canvas={yshift=.5ex}] node {\Bs} (a2)
(a2) edge[transform canvas={yshift=-.4ex}] node {\BIn} (a1)
;
\end{tikzpicture}\]
with $\BIn$ denoting the subcategory inclusion, would-be unit $\righta \cn 1 \to \BIn\Bs$, and would-be counit $\lefta \cn \Bs\BIn \to 1$.  For a bipermutative category $\A$, each of the two composites
\[\begin{split}
\A \fto{\righta} \BsA & \fto{\lefta} \A \andspace \\
\BsA \fto{\Bs\righta} \Bs\BsA & \fto{\lefta} \BsA
\end{split}\]
is the identity symmetric bimonoidal functor by \cref{rl_id} and \cref{def:Bsf,def:lefta,def:righta}.  However, the would-be counit $\lefta \cn \Bs\BIn \to 1$ is \emph{not} well defined because the symmetric bimonoidal functor $\lefta \cn \BsA \to \A$ is not strict.  As explained in \cref{def:lefta}, $\lefta$ is additively strict and multiplicatively strictly unital strong, with the multiplicative monoidal constraint $\leftatwotimes$ given by a Laplaza coherence isomorphism \cref{leftatwotimes}.  Thus, in the bipermutative case, the quadruple $(\Bs,\BIn,\righta,\lefta)$ is not an adjunction.
\end{remark}

\section{Application to Multiplicative Infinite Loop Space Theory}\label{sec:applications}

This section explains that, in the context of multiplicative infinite loop space theory, \cref{thm:sbfstrict} allows us to restrict from multiplicatively strong symmetric bimonoidal functors \cref{mult_strong} between bipermutative categories to \emph{strict} symmetric bimonoidal functors.  The first half of this section recalls May's multiplicative machine.  The application of \cref{thm:sbfstrict} is discussed in the second half of this section.  \cref{rk:Kem} explains why \cref{thm:sbfstrict} does \emph{not} directly apply to Elmendorf-Mandell multifunctorial $K$-theory.

\subsection*{May's Multiplicative Infinite Loop Space Machine}

This subsection briefly recalls May's multiplicative machine $\mayk$ from \cite{may-construction,may-precisely}, which update and correct earlier treatment in \cite{may-einfinity,may-multiplicative}.  For our application below, we only need to know some formal properties of $\mayk$.  The reader is referred to \cite{may-construction,may-precisely} for details.  May's multiplicative machine goes from the category $\biperm$ of small bipermutative categories and symmetric bimonoidal functors \pcref{def:bipermcat,def:sbfunctor} to the category $\Einfrsp$ of $\Einf$-ring spectra.  It is a composite of several functors as follows.
\begin{equation}\label{may-machine}
\begin{tikzpicture}[vcenter]
\def\v{-1.3} \def\h{3.2}
\draw[0cell=1]
(0,0) node (a1) {\biperm}
(a1)++(0,\v) node (a2) {\Fifcat}
(a2)++(\h,0) node (a3) {\Fiftop}
(a3)++(\h,0) node (a4) {\Einfrtop}
(a4)++(0,-\v) node (a5) {\Einfrsp}
;
\draw[1cell=.9]
(a1) edge node {\mayk} (a5)
(a1) edge node[swap] {\mayj} (a2)
(a2) edge node {\cla} (a3)
(a3) edge node {\mayl} (a4)
(a4) edge node[swap] {\maye} (a5)
;
\end{tikzpicture}
\end{equation}
We briefly discuss each functor that comprises $\mayk$.  

\parhead{The functor $\mayj$}.  In the first functor $\mayj$, $\Fsk$ is the category of pointed finite sets of the form $\ordn = \{0,1,\ldots,n\}$ with basepoint $0$ for $n \geq 0$, and pointed functions.  There is a functor 
\[\Fsk \fto{\mayp} \Cat\] 
that sends an object $\ordn$ to the product category $\Fsk^n$.  Recalling that $\ufsm$ is the unpointed finite set $\{1,2,\ldots,m\}$, the value of $\mayp$ at a pointed function $f \cn \ordm \to \ordn$ is the functor
\[\mayp\ordm = \Fsk^m \fto{\mayp f} \Fsk^n = \mayp\ordn\]
that sends an $m$-tuple $\ang{\ordr_i}_{i \in \ufsm} \in \Fsk^m$ of pointed finite sets or pointed functions to the $n$-tuple
\[\bang{\sma_{i \in \finv(j)} \ordr_i}_{j \in \ufsn} \in \Fsk^n.\]
For pointed finite sets $\ord{a}$ and $\ord{b}$, the smash product $\ord{a} \sma \ord{b}$ is identified with $\ord{ab}$ using the lexicographic ordering \cref{lex_order}.  An empty smash product is $\ord{1}$.  The category $\FiF$ is the wreath product of the functor $\mayp$ \cite[Def.\! 5.1]{may-construction}.  

The diagram category $\Fifcat$ has functors $\FiF \to \Cat$ as objects, called $\FiF$-categories, and natural transformations as morphisms.  The first functor $\mayj$ in \cref{may-machine}, discussed in \cite[Section 13]{may-construction}, is a composite of two functors.
\begin{itemize}
\item The first part of $\mayj$ sends small bipermutative categories and symmetric bimonoidal functors to colax pseudofunctors $\FiF \to \Cat$ and lax transformations \cite[4.1.2 and 4.2.1]{johnson-yau}. 
\item The second part of $\mayj$ applies Street's strictification \cite{street_lax}, in the colax form in \cite[3.4]{may-pairings}, to strictify colax pseudofunctors and lax transformations to $\FiF$-categories and natural transformations.
\end{itemize}

\parhead{The functor $\cla$}.  The second functor 
\[\Fifcat \fto{\cla} \Fiftop\] 
in \cref{may-machine} is given levelwise by the classifying space functor $\cla \cn \Cat \to \Top$, from the category $\Cat$ of small categories and functors to the category $\Top$ of compactly generated weak Hausdorff spaces and continuous morphisms.

\parhead{The functor $\mayl$}.  The third functor 
\[\Fiftop \fto{\mayl} \Einfrtop\] 
in \cref{may-machine}, discussed in \cite[Sections 1--10]{may-construction}, sends functors $\FiF \to \Top$ and natural transformations to $\Einf$-ring spaces and their morphisms.  An \emph{$\Einf$-operad pair} $(\mayc,\mayg)$ consists of two topological $\Einf$-operads, $\mayc$ and $\mayg$, together with a $\mayg$-action on $\mayc$ of the form
\[\mayg(m) \times \prod_{i \in \ufsm} \mayc(r_i) \fto{\la} \mayc(r_1 \cdots r_m)\]
that satisfies some associativity, unity, equivariance, and distributivity axioms \cite[Def.\! 4.2]{may-construction}.  The $\Einf$-operad $\mayc$ parametrizes addition, and the $\Einf$-operad $\mayg$ parametrizes multiplication.  As discussed in \cite[Sections 2 and 3]{may-precisely}, the canonical choices for $\mayg$ and $\mayc$ are, respectively, the linear isometries operad and the infinite Steiner operad \cite{steiner}.  Given an $\Einf$-operad pair $(\mayc,\mayg)$, an \emph{$\Einf$-ring space}, as defined in \cite[Section 1]{may-precisely}, is a space equipped with a $\mayc$-action and a $\mayg$-action, and it satisfies a distributivity axiom.  

The functor $\mayl$ is a composite of three functors: a pullback functor followed by two two-sided monadic bar constructions.  Each of the three constituent functors of $\mayl$ induces an equivalence of homotopy categories.  See \cite[3.10, 5.11, 8.6, and 10.6]{may-construction}.

\parhead{The functor $\maye$}.  The last functor 
\[\Einfrtop \fto{\maye} \Einfrsp\]
in \cref{may-machine} sends $\Einf$-ring spaces to $\Einf$-ring spectra \cite[Section 5]{may-precisely}.  It is given by another two-sided monadic bar construction \cite[9.12]{may-precisely}.  In fact, $\maye$ is a specialization of the additive infinite loop space machine in \cite[9.3]{may-precisely}.  It induces an equivalence between the homotopy category of ringlike $\Einf$-ring spaces and the homotopy category of connective $\Einf$-ring spectra.  We note that $\Einf$-ring spaces and $\Einf$-ring spectra are \emph{not} defined as algebras over some operads in some symmetric monoidal categories.

In summary, at the object level, May's multiplicative infinite loop space machine $\mayk$ sends each small bipermutative category to an $\FiF$-category (by $\mayj$), then an $\FiF$-space (by $\cla$), then an $\Einf$-ring space (by $\mayl$), and finally an $\Einf$-ring spectrum (by $\maye$).

\subsection*{Application of \cref{thm:sbfstrict}}

Now we discuss how \cref{thm:sbfstrict} applies to May's multiplicative machine $\mayk = \maye\mayjbl$ \cref{may-machine}.  Suppose $\sff \cn \A \to \B$ is a multiplicatively strong symmetric bimonoidal functor \cref{mult_strong} between small bipermutative categories.  By \cref{thm:sbfstrict} \eqref{thm:sbfstrict_iii} and \cref{righta_nat_diagram}, there is a commutative diagram of multiplicatively strong symmetric bimonoidal functors as follows.
\begin{equation}\label{righta_nat}
\begin{tikzpicture}[vcenter]
\def\v{-1.3}
\draw[0cell]
(0,0) node (a11) {\A}
(a11)++(1.8,0) node (a12) {\BsA}
(a11)++(0,\v) node (a21) {\B}
(a12)++(0,\v) node (a22) {\BsB}
;
\draw[1cell=.9]
(a11) edge node {\righta} (a12)
(a12) edge node {\Bsf} (a22)
(a11) edge node[swap] {\sff} (a21)
(a21) edge node {\righta} (a22)
;
\end{tikzpicture}
\end{equation}
The multiplicatively strong symmetric bimonoidal functor $\righta$ is constructed in \cref{def:righta}, and the strict symmetric bimonoidal functor $\Bsf$ is constructed in \cref{def:Bsf}.  Our claim is this:
\begin{quote}
In the context of May's multiplicative infinite loop space machine, $\sff$ can be replaced by $\Bsf$ via the diagram \cref{righta_nat}.
\end{quote}

In more details, applying the composite functor $\mayjbl$ to the diagram \cref{righta_nat} yields the commutative diagram 
\begin{equation}\label{lbjr}
\begin{tikzpicture}[vcenter]
\def\v{-1.3}
\draw[0cell]
(0,0) node (a11) {\mayjbl\A}
(a11)++(3.2,0) node (a12) {\mayjbl\BsA}
(a11)++(0,\v) node (a21) {\mayjbl\B}
(a12)++(0,\v) node (a22) {\mayjbl\BsB}
;
\draw[1cell=.9]
(a11) edge node {\mayjbl\righta} node[swap] {\sim} (a12)
(a12) edge[transform canvas={xshift=-1em}] node {\mayjbl\Bsf} (a22)
(a11) edge[transform canvas={xshift=1ex}] node[swap] {\mayjbl\sff} (a21)
(a21) edge node {\mayjbl\righta} node[swap] {\sim} (a22)
;
\end{tikzpicture}
\end{equation}
in $\Einfrtop$.  By \cref{thm:sbfstrict} \eqref{thm:sbfstrict_ii}, $\righta$ has a left adjoint $\lefta \cn \BsA \to \A$ that is an additively strict and multiplicatively strictly unital strong symmetric bimonoidal functor, and the same holds for $\B$.  Due to the presence of the levelwise classifying space functor $\cla$, each of the two $\Einf$-ring space morphisms $\mayjbl\righta$ in the diagram \cref{lbjr} has an underlying weak homotopy equivalence.  Thus, $\mayjbl\sff$ can be replaced by $\mayjbl\Bsf$.

Further applying the functor $\maye$ to the diagram \cref{lbjr} yields the commutative diagram
\begin{equation}\label{maykr}
\begin{tikzpicture}[vcenter]
\def\v{-1.3}
\draw[0cell]
(0,0) node (a11) {\mayk\A}
(a11)++(2.3,0) node (a12) {\mayk\BsA}
(a11)++(0,\v) node (a21) {\mayk\B}
(a12)++(0,\v) node (a22) {\mayk\BsB}
;
\draw[1cell=.9]
(a11) edge node {\mayk\righta} node[swap] {\sim} (a12)
(a12) edge node {\mayk\Bsf} (a22)
(a11) edge node[swap] {\mayk\sff} (a21)
(a21) edge node {\mayk\righta} node[swap] {\sim} (a22)
;
\end{tikzpicture}
\end{equation}
 in $\Einfrsp$.  By the homotopical properties of $\maye$, each of the two $\Einf$-ring spectrum morphisms $\mayk\righta$ in the diagram \cref{maykr} has an underlying weak equivalence.  Thus, $\mayk\sff$ can be replaced by $\mayk\Bsf$.

In summary, using May's multiplicative infinite loop space machine from small bipermutative categories to either $\Einf$-ring spaces or $\Einf$-ring spectra, there is no loss of generality in replacing multiplicatively strong symmetric bimonoidal functors $\sff$ by \emph{strict} symmetric bimonoidal functors $\Bsf$ via the diagram \cref{righta_nat}.

We end with the following remark about another multiplicative infinite loop space machine.

\begin{remark}[Elmendorf-Mandell $K$-Theory]\label{rk:Kem}
The application of \cref{thm:sbfstrict}, as discussed in \cref{righta_nat,lbjr,maykr}, applies to May's multiplicative infinite loop space machine.  That discussion cannot be naively replicated using Elmendorf-Mandell $K$-theory \cite{elmendorf-mandell}, denoted $\Kem$.  The latter is an \emph{enriched multifunctor} from the categorically enriched multicategory of small permutative categories to the simplicially enriched multicategory of symmetric spectra.  By multifunctoriality, $\Kem$ transports small bipermutative categories to $\Einf$-algebras in symmetric spectra \cite[Cor.\ 3.9]{elmendorf-mandell}; see \cite[Section 11.6]{cerberusIII} for a detailed exposition.  We note that Elmendorf and Mandell's definition of a bipermutative category \cite[Def.\ 3.6]{elmendorf-mandell} is more general than May's in \cref{def:bipermcat}.  Moreover, $\Einf$-algebras in symmetric spectra are not defined in the same way as May's $\Einf$-ring spectra \cite[Def.\ 5.3]{may-precisely}.  At the morphism level, $\Kem$ takes \emph{strictly unital} symmetric monoidal functors as input.  Thus, even if $\sff$ is additively strictly unital, $\Kem$ does not directly apply to $\righta$ and the diagram \cref{righta_nat}, since the additive unit constraint $\rightazeroplus$ \cref{rightazeroplus} is not even an isomorphism.  

On the bright side, it is conceivable that there is an equivalence of homotopy theories
\[\begin{tikzpicture}
\def\h{2.75}
\draw[0cell]
(0,0) node (a1) {\biperm}
(a1)++(\h,0) node (a2) {\bipermasu}
;
\draw[1cell=.9]
(a1) edge[transform canvas={yshift=.5ex}] node {\Su} (a2)
(a2) edge[transform canvas={yshift=-.4ex}] node {I} (a1)
;
\end{tikzpicture}\]
where $I$ denotes the inclusion from the subcategory $\bipermasu$ of additively strictly unital symmetric bimonoidal functors.  Using the hypothetical functor $\Su$, the composite $\Kem\Su$ can then be applied to the diagram \cref{righta_nat} to possibly yield a diagram of $\Einf$-algebras in symmetric spectra, similar to \cref{maykr}.  For permutative categories, an equivalence of homotopy theories analogous to $(\Su,I)$ is constructed in \cite[Def.\ 9.6 and Prop.\ 10.2]{johnson-yau-Fmulti}.
\end{remark}

\bibliographystyle{sty/amsalpha3}
\bibliography{references}

\end{document}